\documentclass[preprint, 12pt,a4paper,reqno]{elsarticle}
\usepackage[utf8]{inputenc}
\usepackage[english]{babel}
\usepackage{hyperref}

\usepackage{pdflscape}

\usepackage{amsmath, amsthm, amsfonts}
\usepackage{multirow} 
\usepackage{lineno}
\usepackage{booktabs}
\usepackage{graphicx}
\usepackage{array}
\usepackage{longtable}
\usepackage{tikz-cd}
\usepackage{url}
\usepackage{hyperref}

\usepackage{amsmath}
\usepackage{amssymb, amscd}
\usepackage[all, cmtip]{xy}
\usepackage{xcolor}
\usepackage{multirow} 
\usepackage{longtable}
\usepackage{array}

\newproof{pf}{Proof}
\newproof{pot}{Proof of Theorem \ref{thm2}}

\newtheorem{thm}{Theorem}[section]
\newtheorem{conjecture}[thm]{Conjecture}

\newtheorem{lemma}[thm]{Lemma}

{\sc}
\newtheorem{example}{Example}
{}
{}
{\sc}
\newdefinition{rmk}{Remark}[section]
\newdefinition{remarks}{Remarks}[section]

\newcommand{\thismonth}{\ifcase\month\or
  January\or February\or March\or April\or May\or June\or
  July\or August\or September\or October\or November\or December\fi
  \space\number\year}

\makeatletter
\newcommand{\rssymb}[2]{\newcommand{#1}{{\mathrmsl{#2}}}}
\newcommand{\calsymb}[2]{\newcommand{#1}{{\mathcal{#2}}}}
\newcommand{\bbsymb}[2]{\newcommand{#1}{{\mathbb{#2}}}}
\newcommand{\lieoper}[2]{\newcommand{#1}{\mathop
  {\mathfrak{#2}\null}\nolimits}}
\newcommand{\oper}[3][n]{\newcommand{#2}{\mathop
  {\mathrm{#3}\null}\ifx n#1\nolimits\else\limits\fi}}
\newcommand{\rsoper}[3][n]{\newcommand{#2}{\mathop
  {\mathrmsl{#3}\null}\ifx n#1\nolimits\else\limits\fi}}
\bbsymb\C{C} \bbsymb\F{F} \bbsymb\HQ{H}\bbsymb\N{N} \bbsymb\Q{Q}
\bbsymb\R{R} \bbsymb\U{U} \bbsymb\V{V} \bbsymb\W{W} \bbsymb\Z{Z}
\bbsymb\bbf{F} \bbsymb\bbk{K} \bbsymb\bbi{I} \bbsymb\bbl{L}
\bbsymb\bbo{O} \bbsymb\bbj{J} \bbsymb\bby{Y} \bbsymb\bbp{P}
\bbsymb\bba{A}
\calsymb\cA{A} \calsymb\cB{B} \calsymb\cC{C} 
\calsymb\cM{M} \calsymb\cN{N} \calsymb\cO{O} \calsymb\cP{P}
\calsymb\cU{U} \calsymb\cV{V} \calsymb\cW{W} \calsymb\cX{X}
\calsymb\cY{Y} \calsymb\cZ{Z}
\renewcommand{\geq}{\geqslant} \renewcommand{\leq}{\leqslant}
\oper\End{End} \oper\Hom{Hom}                    
\oper\Sym{Sym} \oper\Skew{Skew}
\oper\Aut{Aut}                                   
\oper\GL{GL} \oper\SL{SL}\oper\Symp{Sp} \oper\CO{CO} \oper\On{O}
\oper\SO{SO} \oper\Pin{Pin} \oper\Spin{Spin} \oper\CU{CU}
\oper\Un{U} \oper\SU{SU} \oper\PSU{PSU} \rsoper\Diff{Diff}
\rsoper\SDiff{SDiff}
\lieoper\der{der}                                
\lieoper\gl{gl} \lieoper\sgl{sl}\lieoper\symp{sp} \lieoper\co{co}
\lieoper\so{so} \lieoper\spin{spin} \lieoper\cu{cu} \lieoper\un{u}
\lieoper\su{su} \rsoper\Vect{Vect} \rsoper\Ham{Ham}
\def\la#1{\hbox to #1pc{\leftarrowfill}}
\def\ra#1{\hbox to #1pc{\rightarrowfill}}

\newcommand{\Norm}[2][]{\bigl|\mkern-3mu\bigr|#2\bigr|\mkern-3mu\bigr|
  _{\lower1pt\hbox{${}_{#1}$}}}

\rsoper\dimn{dim}                           
\rsoper\grad{grad}                          
\rsoper\kernel{ker}\rsoper\image{im}        
\rsoper\alt{alt}   \rsoper\sym{sym}         
\rsoper\Ad{Ad}     \rsoper\ad{ad}           
\rsoper\CoAd{CoAd} \rsoper\coad{coad}       
\rsoper\trace{tr}  \rsoper\trfree{tf}       
\rsoper\detm{det}                           
\rsoper\Vol{Vol}                            
\rsoper\divg{div}                           
\rsoper\sign{sign}                          
\rssymb\iden{id}                            
\rssymb\vol{vol}                            
\oper\Imag{Im}\oper\Real{Re}                
\newcommand{\sd}{{\raise1pt\hbox{$\scriptscriptstyle +$}}}
\newcommand{\asd}{{\raise1pt\hbox{$\scriptscriptstyle -$}}}
\newcommand{\sdasd}{{\raise1pt\hbox{$\scriptscriptstyle\pm$}}}

\newcommand{\asdsd}{{\raise1pt\hbox{$\scriptscriptstyle\mp$}}}

\rsoper\scal{scal}
\def\kahl/{k\"ahler}
\def\Kahl/{K{\"a}hler}

\def\@author#1{\g@addto@macro\elsauthors{\normalsize%
    \def\baselinestretch{1}%
    \upshape\authorsep#1\unskip\textsuperscript{%
      \ifx\@fnmark\@empty\else\unskip\sep\@fnmark\let\sep=,\fi
      \ifx\@corref\@empty\else\unskip\sep\@corref\let\sep=,\fi
      }%
    \def\authorsep{\unskip,\space}%
    \global\let\@fnmark\@empty
    \global\let\@corref\@empty  
    \global\let\sep\@empty}%
    \@eadauthor={#1}
}

\begin{document}

\begin{frontmatter}

\title{Non-existence of Extremal Sasaki metrics and the Berglund-H\"ubsch transpose}

\author[add1]{Jaime Cuadros Valle\corref{cor1}%
}
\ead{jcuadros@pucp.edu.pe}
\address[add1]{Departamento de Ciencias, Sección Matemáticas,
 Pontificia Universidad Católica del Perú, Av. Universitaria 1801, Lima 32, Perú}

\author[add2]{Ralph R. Gomez}
\ead{rgomez1@swarthmore.edu}
\address[add2]{Department of Mathematics and Statistics, Swarthmore College,
500 College Avenue, Swarthmore,
Swarthmore, PA 19081, USA}

\author[add1]{Joe Lope Vicente}
\ead{j.lope@pucp.edu.pe}

\cortext[cor1]{Corresponding author}

\date{\thismonth}

\begin{abstract}

We use the Berglund-Hübsch transpose rule from classical mirror symmetry in the context of Sasakian geometry \cite{CGL} and results on relative K-stability in the Sasaki setting developed by Boyer and van Coevering in  \cite{BvC} to exhibit  examples  of Sasaki manifolds with big Sasaki cones that have no extremal Sasaki metrics at all.   Previously, examples with this feature were  produced in  \cite{BvC}  for Brieskorn-Pham polynomials or their deformations. Our examples are based on the more general framework of invertible polynomials. In particular, we  construct families  of links that preserve the emptiness of the extremal Sasaki-Reeb cone via the Berglund-Hübsch rule: if the link does not admit extremal Sasaki metrics then its  Berglund-H\"ubsch dual preserves this property and moreover this dual admits a representative in its local moduli with a larger Sasaki-Reeb cone which remains  obstructed to admitting  extremal Sasaki metrics. Some of the examples exhibited here  have the homotopy type of a sphere or are rational homology spheres. 
\end{abstract}


\medskip

\medskip

\vspace{-2mm} 

\end{frontmatter}


\section{Introduction} 
Introduced in \cite{BGS}, extremal Sasaki metrics are critical points  of a modified version of the  Calabi functional of K\"ahler geometry \cite{Ca}, this time defined over a suitable space of Sasaki metrics which fixes the Reeb vector field. These critical points turn out to be Sasaki metrics whose transversal scalar curvatures are such that their gradients are transversely holomorphic. In particular, Sasaki metrics with constant scalar curvature and hence 
Sasaki-Einstein metrics are extremal. Thus, these metrics  provide  an appropriate  notion of  canonical metrics. 

In order to determine the set of  extremal Sasaki metrics, called the extremal cone,  one considers the Sasaki-Reeb cone of a Sasakian manifold, which consists of  the deformations of  the given Reeb vector field by Hamiltonian holomorphic vector fields. 
It is natural to ask whether  the set of extremal Sasaki metrics is empty or not, and if not  when it consists of the whole Sasaki-Reeb cone or if among them there exist metrics with constant scalar curvature or even Sasaki-Einstein metrics. 
The first examples of Sasaki manifolds with Sasaki-Reeb cone of dimension larger than one, where the extremal set is empty were exhibited in \cite{BvC}. These examples were given as links of Brieskorn-Pham polynomials and as such, their Sasaki-Reeb cones are of Gorenstein type, that is, the polarization is induced  by the canonical or anticanonical bundle. In \cite{BHLTF}, Boyer et al.  give examples with the same feature as in \cite{BvC} for  the non-Gorenstein case; moreover the authors found examples of Sasaki manifolds admitting extremal metrics  with no representative admitting constant scalar curvature. 

In this paper, we continue studying the Gorenstein case.  The novelty in our heuristic approach is that the examples presented here  are derived from links of the more general framework of invertible polynomials, that is, weighted homogeneous polynomials $f \in \mathbb{C}\left[x_0, \ldots, x_n\right]$, which are a sum of exactly $n+1$ monomials, such that the weights $w_0, \ldots, w_n$ of the variables $x_0, \ldots, x_n$ are unique  and that  its affine cone is smooth outside  $(0, \ldots, 0),$ {\it i.e.,} $f$ is quasismooth. 
Brieskorn-Pham polynomials are an important building block for invertible polynomials.  This type of polynomials has been intensely studied in the context of Sasakian geometry and  as a result it is known that  homotopy spheres given as links of Brieskorn-Pham polynomials admit infinitely many Sasaki-Einstein metrics \cite{GK,BGK, ST,LST}. The focus on Brieskorn-Pham polynomials were motivated by the well-known fact that every homotopy sphere that bounds a parallelizable manifold is the link of some Brieskorn-Pham polynomial in the sense that every differentiable structure on these spheres can be obtained as a link of a this type of singularity \cite{Br}. However, from a Riemannian viewpoint it could be of benefit to obtain this differentiable structure on a homotopy sphere from a link of an invertible polynomial that cannot be written as a Brieskorn-Pham polynomial and try to determine whether this link admits a extremal Sasaki metrics. So far, efforts to establish the existence of extremal Sasaki  metrics in this kind of links have not been successful. We conjecture that if a link $L_f$ is such  that the following conditions are satisfied 
\begin{itemize}
\item $L_f$ is given by an invertible polynomial that cannot be written as a Brieskorn-Pham polynomial.  
\item $L_f$ is a  homotopy sphere 
\end{itemize}
then $L_f$ does not admit extremal Sasaki metrics in its whole Sasaki-Reeb cone. Theorem 4.3(3)  provides some mild evidence for the plausability of this statement. Actually some of the examples we produce are exotic spheres (see Table 1).

Our strategy allows us to make use of the Berglund-Hübsch transpose rule of classical mirror symmetry \cite{BH}, whose first use in the context of Sasakian geometry was done in \cite{Go}.  Recall that the Berglund-Hübsch rule considers an invertible polynomial $f=\sum_{i=0}^n \prod_{j=0}^n x_j^{a_{i j}}$ cutting out an orbifold of degree $d$ in $\mathbb{P}(\bf w)$ and defines the transpose polynomial $f^T$  by transposing the exponential  matrix $A=\left(a_{i j}\right)_{i, j}$ of the original polynomial, that is,  $f^T=\sum_{i=0}^n \prod_{j=0}^n x_j^{a_{j i}},$ a polynomial that cuts out an orbifold of degree $\tilde{d}$ in $\mathbb{P}(\tilde{\bf w}).$ Then one considers the links 
$L_f(\mathbf{w}, d)$ and  $L_{f^{T}}(\tilde{\mathbf{w}}, \tilde{d})$ associated to each of these polynomials. We will sometimes say that these two links are {\it Berglund-H\"ubsch duals} to one another. Note that for Brieskorn-Pham  polynomials the Berglund-H\"ubsch rule does not produce new data since their exponential matrices are given by  diagonal matrices with integer coefficients so this type of polynomials are self-dual. This method has proven fruitful in the context of Sasakian geometry, for instance in \cite{CGL} we  found conditions for rational homology 7-spheres to preserve the existence of extremal Sasaki metrics under the Berglund-H\"ubsch rule, however in those cases the dimension of the Sasaki-Reeb cone remained constant and of dimension one.
The following diagram succinctly summarizes the procedure described above (see details in Section 3), where BH denotes the Berglund-H\"ubsch transpose rule:

\begin{center}
\begin{tikzcd}
{f=0} \arrow{d}\arrow{r}{BH} & {f^T=0}\arrow{d}\\
L_f(\mathbf{w}, d)  \arrow{r}{BH} & L_{f^{T}}(\tilde{\mathbf{w}}, \tilde{d}).\\
\end{tikzcd}
\end{center}

In this work we also study whether the  Berglund-H\"ubsch transpose preserves the emptiness of the extremal Sasaki-Reeb cone. In Subsections  4.2 and 4.3 we exhibit families of links $L_f(\mathbf{w}, d)$ with  Sasaki-Reeb cones of dimension 1 and no extremal Sasaki metrics on them for which the  Berglund-H\"ubsch rule produces links with Sasaki-Reeb cones of dimension greater than one and with empty extremal Sasaki-Reeb cone. More precisely, the dual link $L_{f^{T}}(\tilde{\mathbf{w}}, \tilde{d})$ admits a {\it perturbation} in its local moduli, that is, an element in the complex space of weighted homogeneous polynomials of degree $\tilde{d}$,  $H^0(\mathbb{P}(\tilde{\mathbf{w}}),\mathcal{O}(\tilde{d})),$ with  a larger group of automorphism which  produces a link with larger Sasaki-Reeb cone that remains obstructed to admitting  extremal Sasaki metrics.

In order to  achieve the results given in this paper, we apply results of  Boyer and van Coevering  in \cite{BvC}, where the authors used a relative notion of  K-stability of Collins and Székelyhidi \cite{CS} to generalize the Lichnerowicz inequality of Gauntlett, Martelli, Sparks and Yau\cite{GMSY}, inequality that gives an obstruction to the existence of Sasaki-Einstein metrics on links. This generalized Lichnerowicz inequality gives an obstruction to the existence of extremal Sasaki metrics on certain links.

The paper is organized as follows: in Section 2 we present some background material which includes some basics on Sasaki geometry and focus on  links of hypersurface singularities as canonical examples; we also review  how to calculate the homology of the links. In Section 3 we briefly describe the Berglund-Hübsch transpose rule from mirror symmetry and its application to Sasakian geometry. In Section 4 we state and prove our main results and exhibit the examples described above.

\section{Preliminaries} 

\subsection{Sasaki-Reeb cone and extremal Sasaki metrics}
Let us begin by recalling the definition of a Sasakian structure (see \cite{BBG} for a comprehensive survey on this topic). Sasakian geometry is a special type of contact metric structure on a $(2n+1)$-dimensional manifold $M$ described by the quadruple of tensors $(\Phi, \xi, \eta, g)$  such that $\eta$ is a contact 1-form, $\Phi$ is an endomorphism of the tangent bundle, 
$g=d\eta \circ (\mathbb{I}\times \Phi)+\eta\otimes \eta$ is a Riemannian metric and $\xi$ is the Reeb vector field which is Killing. Moreover, the underlying CR-structure $\left(\mathcal{D},\left.\Phi\right|_{\mathcal{D}}\right)$ is integrable, where $\mathcal{D}=\ker \eta.$

A Sasaki manifold has a transverse  Kählerian structure  $(\nu(\mathcal{F}_\xi), \bar{J})$ 
determined by  the normal bundle $\nu(\mathcal{F}_\xi)$ of the characteristic foliation $\mathcal{F}_\xi$ determined by $\xi$ and the natural complex structure $\bar{J}$ in $\nu(\mathcal{F}_\xi)$ induced by $\left.\Phi\right|_{\mathcal{D}}.$ In fact,  when all the orbits of $\xi$ are closed, the Reeb vector field $\xi$ generates a locally free circle action whose quotient is a Kähler  orbifold. In this case the Sasakian structure is called quasi-regular. When the action is free it is called  regular and the quotient is a Kähler manifold. In the irregular case when the orbits of the Reeb vector field  are not closed the local quotients are Kähler. 

For the next definitions we follow \cite{BGS} and \cite{BHLTF}. Let $\mathcal{S}=(\xi, \eta, \Phi, g)$ be a $T$-invariant Sasakian structure where $T$ is a compact torus acting effectively on $M$, that is, $T \subset \operatorname{Diff}(M)$  such that $\xi \in \operatorname{Lie} T=: \mathfrak{t}$. The Sasaki-Reeb cone of $(M, \mathcal{S})$, relative to $T$, is given by  
$$
\mathfrak{t}^{+}(\mathcal{D}, J):=\{\tilde{\xi} \in \mathfrak{t} \mid \eta(\tilde{\xi})>0\}.
$$
Each $\tilde{\xi} \in \mathfrak{t}^{+}(\mathcal{D}, J)$ determines a $T$-invariant Sasakian structure, compatible with the CR-structure $(\mathcal{D}, J)$ with contact form $\displaystyle{\frac{\eta}{\eta(\tilde{\xi})}}$. Thus, elements of $\mathfrak{t}^{+}(\mathcal{D}, J)$ are called Sasakian metrics. We will   
consider $\mathfrak{t}^{+}(\mathcal{D}, J)$  as a $k$-dimensional smooth family of Sasakian structures. It is often convenient to use the complexity $n+1-k$ instead of  $\operatorname{dim} \mathfrak{t}^{+}=k$. Then the complexity is an integer from 0 to $n$ with 0 being the toric case.

Let us denote the space of Sasakian structures with Reeb vector field $\xi$ and transverse holomorphic structure $\bar{J}$ by $\mathcal{S}(\xi, \bar{J})$, that is, 
the set  with the same complex normal bundle $(\nu(\mathcal{F}_\xi), \bar{J})$. This is an infinite dimensional Fréchet space, so it is convenient to consider  the space of  $T$-invariant isotopic contact structures $\mathcal{S}(\xi, \bar{J})^T.$  For any 
$(\eta, \xi, \Phi, g) \in \mathcal{S}(\xi, \bar{J})$ on a manifold $M,$ the Calabi functional \cite{BGS} is defined by 
$$
E(g)=\int_M s_g^2 d \mu_g,
$$
where $s_g$ is the scalar curvature. A Sasaki metric $g$ is extremal if it is a critical point of $E$. Equivalently, the metric $g$ is extremal if the $(1,0)$ component of the gradient vector field $\partial^\# s_g$ is a transversely holomorphic vector field. Clearly, this  notion is a generalization of Sasaki metrics with constant scalar curvature and hence of Sasaki-Einstein metrics. One defines the extremal Sasaki-Reeb cone as the set 
{\small{
$$\mathfrak{e}(\mathcal{D}, J):=\left\{\xi \in \mathfrak{t}^{+}(\mathcal{D}, J) \mid  \hbox{ there is an extremal Sasakian structure in } \mathcal{S}(\xi, \bar{J})^T\right\}.$$}}

Explicit constructions of Sasakian manifolds in a very general way can be given as total spaces of $S^1$-orbibundles over projective orbifolds. As a special case, one can build up Sasakian structures on  links of  hypersurface singularities of weighted homogeoneous polynomials and we explain how to do this below.   

\subsection{Links of hypersurface singularities and Sasakian geometry} 
Let us recall  that a polynomial $f \in \mathbb{C}\left[z_0, \ldots, z_n\right]$ is said to be a weighted homogeneous polynomial of degree $d$ and weight vector $\mathbf{w}=$ $\left(w_0, \ldots, w_n\right)$, with weights $w_i \in \mathbb{Z}^{+}$, and if for any $\lambda \in \mathbb{C}^{\times}$
$$
f\left(\lambda^{w_0} z_0, \ldots, \lambda^{w_n} z_n\right)=\lambda^d f\left(z_0, \ldots, z_n\right).
$$

We are interested in those weighted homogeneous polynomials $f$ whose zero locus in $\mathbb{C}^{n+1}$ has only an isolated singularity at the origin. A link, denoted by $L_f(\mathbf{w}, d)$ of a hypersurface singularity is defined  as $f^{-1}(0) \cap S^{2n+1}$, where $S^{2n+1}$ is the $(2n+1)$-sphere in Euclidean space. By the Milnor fibration theorem \cite{Mi}, $L_f(\mathbf{w}, d)$ is a closed ($n-2)$-connected manifold that bounds a parallelizable manifold with the homotopy type of a bouquet of $n$-spheres. Furthermore, $L_f(\mathbf{w}, d)$ admits a quasi-regular Sasaki structure $\mathcal{S}=\left(\xi_{\mathbf{w}}, \eta_{\mathbf{w}}, \Phi_{\mathbf{w}}, g_{\mathbf{w}}\right)$  which is the restriction of the weighted Sasakian structure on the sphere $S^{2 n+1}$ with Reeb vector field $$\xi_{\mathrm{w}}=\sum_k w_k H_k,$$ where 
$$H_k=-i\left(z_k \partial_{z_k}-\bar{z}_k \partial_{\bar{z}_k}\right)=\left ( y_k \partial_{{x}_k}-x_k\partial_{{y}_k} \right ),$$
for $k=0, \ldots n.$

This structure is quasi-regular and its bundle $(\mathcal{D}, J)$ has $c_1(\mathcal{D})=0$. This latter property implies that
$$
c_1\left(\mathcal{F}_{\xi_{\mathbf{w}}}\right)=a\left[d \eta_{\mathbf{w}}\right]_B
$$ 
for some constant $a$, where $\mathcal{F}_{\xi_{\mathrm{w}}}$ is the characteristic foliation. The sign of $a$ determines the negative, null, and positive cases that we shall refer to below which shows that it is an easy matter to obtain Sasaki-Reeb cones of dimension one.

\begin{thm}{\cite{BGS}}
Let $L_f(\mathbf{w}, d)$ be either a negative or null link of an isolated hypersurface singularity with underlying CR-structure $(\mathcal{D}, J)$. Then its Sasaki-Reeb cone $\mathfrak{t}^{+}(\mathcal{D}, J)$ is one-dimensional. 
\end{thm}

It is well-known that the case of positive links is more complicated. However, for the $S^1$-Seifert bundle $\pi:L_f\rightarrow X_f$ the Sasaki-Reeb cone $\mathfrak{t}^{+}(\mathcal{D}, J)$ associated to the natural Sasakian structure on $L_f$ is generated by the Reeb vector field and the lifting of the vector fields in $\mathfrak{aut}(X_f),$ the Lie algebra of $\mathfrak{Aut}(X_f).$ We have the following result whose proof can be found in \cite{BBG}  Lemma 5.5.3.

\begin{lemma} Let $|\mathbf{w}|$ denote the sum of the weights $w_i$'s. 
If $2 w_i<d$ for all but at most one $i$ then the group $\mathfrak{A u t}\left(X_f\right)$ is discrete. If in addition $n \geq 4$ or $n=3$ and $|\mathbf{w}| \neq d$, then the group $\mathfrak{A u t}\left(X_f\right)$ is finite. In particular, the corresponding link $L_f(\mathbf{w}, d)$ of the isolated hypersurface singularity $f=0$ with underlying $\mathrm{CR}$ structure  $(\mathcal{D}, J)$ has a  Sasaki-Reeb cone $\mathfrak{t}^{+}(\mathcal{D}, J)$ of dimension one.
 \end{lemma}

\subsection{Topology of links}
Recall that 
the Alexander  polynomial $\Delta_f(t)$ in \cite{Mi} associated to a link $L_f$ of dimension $2 n-1$ is the characteristic polynomial of the monodromy map $$h_*: H_{n}(F, \mathbb{Z}) \rightarrow H_{n}(F, \mathbb{Z})$$  induced by the circle action on the Milnor fibre $F$. Then, 
$$\Delta_f(t)=\operatorname{det}\left(t {\mathbb I}-h_*\right).$$ 
Now both $F$ and its closure $\bar{F}$ are homotopy equivalent to a bouquet of $n$-spheres $S^n \vee \cdots \vee S^n,$ and the boundary of $\bar{F}$ is the link $L_f$, which is $(n-2)$-connected. The Betti numbers $b_{n-1}\left(L_f\right)=b_{n}\left(L_f\right)$ equal the number of factors of $(t-1)$ in $\Delta_f(t)$. 
$L_f$ is a rational homology sphere if and only if $\Delta_f(1) \neq 0$ and the order of $H_{n-1}\left(L_f, \mathbb{Z}\right)$ equals $|\Delta_f(1)|$.
$L_f$  is homotopy sphere if and only if $\Delta_f(1)= \pm 1.$ 

If  $L_f$ is homeomorphic to $S^{2 n-1}$ one can can use the following criteria to  characterize the smooth structure 
\begin{enumerate}
\item  When $n$ is odd, i.e. $\operatorname{dim} L_f \equiv 1 \bmod 4$, the smooth structure is completely determined by the Kervaire invariant $c\left(F\right) \in \mathbb{Z} / 2 \mathbb{Z}$ (see \cite{Mi} for the precise definition of this invariant).
\item  When $n$ is even, i.e. $\operatorname{dim} L_f \equiv 3 \bmod 4$, the smooth structure is completely determined by the signature of the intersection product
$$
s: H_n\left(F, \mathbb{Z}\right) \times H_n\left(F, \mathbb{Z}\right) \rightarrow \mathbb{Z}.
$$
\end{enumerate}
If $n$ is odd, one has a  remarkable result of Levine \cite{Le}: 
$$
c\left(F\right)=\left\{\begin{array}{lll}
0, & \text { if } \Delta_f(-1) \equiv \pm 1 & \bmod 8 \\
1, & \text { if } \Delta_f(-1) \equiv \pm 3 & \bmod 8,
\end{array}\right. 
$$
where $c(F)=1$ corresponds to an exotic sphere.

Thus to determine whether a $(4n+1)$-dimensional link $L_f$ is an exotic sphere amounts to computing $\Delta_f(t)$.

In the case that $f$ is a weighted homogeneous polynomial, Milnor and Orlik \cite{MO} associate to any monic polynomial $f$ with roots $\alpha_1, \ldots, \alpha_k \in \mathbb{C}^*$ its divisor
$$
\operatorname{div} f=\left\langle\alpha_1\right\rangle+\cdots+\left\langle\alpha_k\right\rangle
$$
thought of as an element of the integral ring $\mathbb{Z}\left[\mathbb{C}^*\right].$ Let  $\Lambda_n=\operatorname{div}\left(t^n-1\right).$ Notice that $\Lambda_1=\operatorname{div}\left(t-1\right)$ is the ring identity element and the relation $\Lambda_a \Lambda_b=\operatorname{gcd}(a, b) \Lambda_{\operatorname{lcm}(a, b)}.$
Then the divisor of $\Delta_f(t)$ is given by
\begin{equation}
\operatorname{div} \Delta_f=\prod_{i=0}^n\left(\frac{\Lambda_{u_i}}{v_i}-\Lambda_1\right),
\end{equation}
where the $u_i's$  and $v_i's$ are given in terms of the degree $d$ of $f$ and the weight vector ${\bf w}=(w_0,\ldots w_n)$ by the equations 
\begin{equation}
u_i=\frac{d}{\operatorname{gcd}\left(d, w_i\right)}, \quad v_i=\frac{w_i}{\operatorname{gcd}\left(d, w_i\right)}.
\end{equation}


In particular if the link $L_{f}(\textbf{w}, d)$ is such that the  weight vector $\textbf{w}$ is of the form 
$$\textbf{w}=(w_{0},w_{1},w_{2},w_{3},w_{4})=(m_{2}v_{0},m_{2}v_{1},m_{2}v_{2},m_{3}v_{3},m_{3}v_{4})$$ 
where the $v_i's$ are given as in the previous equation, with 
$\gcd(m_{2},m_{3})=1$ and $m_{2}m_{3}=d$, one obtains $u_0=u_1=u_2=m_3$ and 
$u_3=u_4=m_2.$ Thus, applying the relations $\Lambda_a \Lambda_b=\operatorname{gcd}(a, b) \Lambda_{\operatorname{lcm}(a, b)}$ in Equation (1), we have  
\begin{equation}
    \operatorname{div}\Delta_f=\alpha({\bf w})\beta({\bf w})\Lambda_d+\beta({\bf w})\Lambda_{m_3}-\alpha({\bf w})\Lambda_{m_2}-\Lambda_1,
    \end{equation}
with $$\alpha(\textbf{w})=\frac{m_{2}}{v_{3}v_{4}}-\frac{1}{v_{3}}-\frac{1}{v_{4}}$$ and 
$$\beta(\textbf{w})=\left (\frac{m_3}{v_0v_1}-\frac{1}{v_1}-\frac{1}{v_0}\right )\left (\frac{m_3}{v_2}- 1\right ) + \frac{1}{v_2}$$ positive  integers depending on the weights. From Equation (3) we have 
 $$\Delta_f(t)=\frac{(t^d-1)^{\alpha({\bf w})\beta({\bf w})}(t^{m_3}-1)^{\beta({\bf w})}}{(t^{m_2}-1)^{\alpha({\bf w})}(t-1)}.$$
Since a rational homology sphere satisfies $\Delta_f(1)\not= 0,$  we have 
 
$$(\alpha({\bf w})+1)(\beta({\bf w})-1)=0,$$ which implies  
\begin{equation}
\beta({\bf w})=\left (\frac{m_3}{v_0v_1}-\frac{1}{v_1}-\frac{1}{v_0}\right )\left (\frac{m_3}{v_2}- 1\right ) + \frac{1}{v_2}=1.
\end{equation}
\medskip

Returning to the general case, it is well-known (see \cite{MO}) that the Milnor number is given by $$\left(\frac{d-w_0}{w_0}\right)\left(\frac{d-w_1}{w_1}\right)\left(\frac{d-w_2}{w_2}\right)\left(\frac{d-w_3}{w_3}\right)\left(\frac{d-w_4}{w_4}\right)=\mu.$$
One also  can calculate the free part of $H_{n-1}\left(L_f, \mathbb{Z}\right)$ from the following formula : 
\begin{equation*}
b_{n-1}\left(L_{f}\right)=\sum(-1)^{n+1-s} \frac{u_{i_{1}} \cdots u_{i_{s}}}{v_{i_{1}} \cdots v_{i_{s}} \operatorname{lcm}\left(u_{i_{1}}, \ldots, u_{i_{s}}\right)},
\end{equation*}
where the sum is taken over all the $2^{n+1}$ subsets $\left\{i_{1}, \ldots, i_{s}\right\}$ of $\{0, \ldots, n\}$. In \cite{Or}, Orlik gave a conjecture  which allows one to determine the torsion of the homology groups  of the link  in terms of the weight of $f.$  

\begin{conjecture}[Orlik]  Let $L_f$ denote a link of an isolated hypersurface singularity defined by a weighted homogenous polynomial $f$ with weight vector $\mathbf{w}=\left(w_0, \ldots w_n\right)$ and degree $d$.
Consider  $\left\{i_{1}, \ldots, i_{s}\right\} \subset\{0,1, \ldots, n\}$ the ordered set of $s$ indices, that is, $i_{1}<i_{2}<\cdots<i_{s}$. Let us denote by I its power set (consisting of all of the $2^{s}$ subsets of the set), and by $J$ the set of all proper subsets. Given a $(2 n+2)$-tuple $(\mathbf{u}, \mathbf{v})=\left(u_{0}, \ldots, u_{n}, v_{0}, \ldots, v_{n}\right)$ of integers given as in Equation (2), let us  define inductively a set of $2^{s}$ positive integers, one for each ordered element of $I$, as follows:
$$
c_{\emptyset}=\operatorname{gcd}\left(u_{0}, \ldots, u_{n}\right),
$$
and if $\left\{i_{1}, \ldots, i_{s}\right\} \subset\{0,1, \ldots, n\}$ is ordered, then
$$
c_{i_{1}, \ldots, i_{s}}=\frac{\operatorname{gcd}\left(u_{0}, \ldots, \hat{u}_{i_{1}}, \ldots, \hat{u}_{i_{s}}, \ldots, u_{n}\right)}{\prod_{J} c_{j_{1}, \ldots j_{t}}} .
$$
Similarly, we also define a set of $2^{s}$ real numbers by
$$
k_{\emptyset}=\epsilon_{n+1},
$$
and
$$
k_{i_{1}, \ldots, i_{s}}=\epsilon_{n-s+1} \sum_{I}(-1)^{s-t} \frac{u_{j_{1}} \cdots u_{j_{t}}}{v_{j_{1}} \cdots v_{j_{t}} \operatorname{lcm}\left(u_{j_{1}}, \ldots, u_{j_{t}}\right)}
$$
where
$$
\epsilon_{n-s+1}= \begin{cases}0 & \text { if } n-s+1 \text { is even } \\ 1 & \text { if } n-s+1 \text { is odd }\end{cases}
$$
respectively. Finally, for any $j$ such that $1 \leq j \leq r=\left\lfloor\max \left\{k_{i_{1}, \ldots, i_{s}}\right\}\right\rfloor$, where $\lfloor x\rfloor$ is the greatest integer less than or equal to $x$, we set
$
d_{j}=\prod_{k_{i_{1}, \ldots, i_{s}} \geq j} c_{i_{1}, \ldots, i_{s}}.$
Then 

\begin{equation*}
H_{n-1}\left(L_{f}, \mathbb{Z}\right)_{\text {tor }}=\mathbb{Z} / d_{1} \oplus \cdots \oplus \mathbb{Z} / d_{r} .
\end{equation*}

\end{conjecture}
\medskip

This conjecture is true  for invertible polynomials \cite{HM}, that is, for iterated Thom-Sebastiani sum of the following  type singularities 
\begin{enumerate}

\item {\it Chain type singularity:}  a  singularity given by a weighted homogeneous polynomial  of the form
$$
f=f\left(x_{1}, \ldots, x_{n}\right)=x_{1}^{a_{1}+1}+\sum_{i=2}^{n} x_{i-1} x_{i}^{a_{i}}
$$
for some $n \in \mathbb{N}$ and some $a_{1}, \ldots, a_{n} \in \mathbb{N}$. 
\item {\it Cycle or loop type singularity:} a  singularity given by a weighted homogeneous polynomial  of the form
$$
f=f\left(x_{1}, \ldots, x_{n}\right)=\sum_{i=1}^{n-1} x_{i}^{a_{i}} x_{i+1}+x_{n}^{a_{n}} x_{1}
$$
for some $n \in \mathbb{Z}_{\geq 2}$ and some $a_{1}, \ldots, a_{n} \in \mathbb{N}$ which satisfy for even $n$ neither $a_{j}=1$ for all even $j$ nor $a_{j}=1$ for all odd $j.$
Notice that  {\it Brieskorn-Pham singularities}, or {\it BP singularities} $$f=f\left(x_{1}, \ldots, x_{n}\right)=\sum_{i=1}^{n} x_{i}^{a_{i}}$$ for some $n \in \mathbb{N}$ and some $a_{1}, \ldots, a_{n} \in \mathbb{Z}_{\geq 2}$ are special cases of Thom-Sebastiani sum of chain type singularities and sometimes we will tacitly use this phrasing.  
\end{enumerate}

\section{Berglund-Hübsch and the Sasaki-Reeb cone}

In general, invertible polynomials can be written as $$f\left(x_1, \ldots, x_n\right)=\sum_{i=1}^n \prod_{j=1}^n x_j^{A_{i j}}$$ with $\operatorname{det}\left(A_{i j}\right) \neq 0$, and so the Berglund-Hübsch transpose rule  \cite{BH}, is given by 
$$
f^T\left(x_1, \ldots, x_n\right)=\sum_{i=1}^n \prod_{j=1}^n x_j^{A_{j i}}, 
$$
that is,  $f^T\left(x_1, \ldots, x_n\right)$ is defined by transposing the exponent matrix $A=\left(A_{i j}\right)_{i, j}$ of the original polynomial.  It is well-known that the Berglund-Hübsch transpose rule of an invertible polynomial remains invertible. It is natural to consider the links associated to $f$ and $f^T$. Let us see how the Berglund-H\"ubsch transpose rule works in the context of Sasakian geometry through an example in dimension five.

\begin{example}
Consider the invertible polynomial of degree $d=57$ given by
$$f=z_{0}^{7}z_{1}+z_{1}^{4}z_{2}+z_{2}^{2}z_{0}+z_{3}^{3}$$
which defines the hypersurface in weighted projective space $\mathbb{P}(7,8,25,19).$ Note that this polynomial is of BP-cycle type. Then we form the matrix of exponents and obtain
\begin{equation*}
A=\begin{bmatrix}
       7 & 1 & 0 & 0\\
       0 & 4& 1 & 0\\
       1 & 0 & 2 & 0\\
       0 & 0 & 0 & 3\\
       \end{bmatrix}.
\end{equation*}
Now, we can take the transpose of this matrix $A^{T}$ and solve $A^{T}\tilde{\mathbf{w}}^{T}=\tilde{D}$
where $\tilde{\mathbf{w}}=(\tilde{w_{0}},\tilde{w_{1}},\tilde{w_{2}}, \tilde{w_{3}})$ is the weight vector associated to $f^{T}$ and $\tilde{D}=(\tilde{d},\tilde{d},\tilde{d},\tilde{d})^{T}$ is the column vector comprised of the degree of $f^{T}.$ It is easy to check that 
$$f^{T}=z_{0}^{7}z_{2}+z_{0}z_{1}^{4}+z_{1}z_{2}^{2}+z_{3}^{3}$$
is of degree $\tilde{d}=57$ which defines a hypersurface in the weighted project space $\mathbb{P}(5,13,22,19).$
The links $L_{f}$ and $L_{f^{T}}$, associated to $f$ and $f^T$ respectively, are known to be Sasaki-Einstein and diffeomorphic to $\#2(S^{2}\times S^{3})$. In fact both of these 5-dimensional Sasaki-Einstein examples appear in Table 1 in \cite{BGN1}. Thus we see that the Berglund-H\"ubsch transpose rule maps a Sasaki-Einstein structure on $2\#(S^{2}\times S^{3})$ to another Sasaki-Einstein structure on that same manifold.
\end{example}

This example motivates the idea that one should perhaps study the Berglund-H\"ubsch transpose rule in the context of the Sasaki-Reeb cone. This type of analysis was conducted more systematically with links with the rational homology of a 7-sphere, where interesting results were found in \cite{CGL}. However, for the cases studied there, the dimension of  Sasaki-Reeb cone of the corresponding links remained fixed and of dimension one. This is not always the case, for instance, certain Thom-Sebastiani sum of two or more  chain type singularities has a Berglund-H\"ubsch dual which admits a deformation or perturbation in its   moduli, a link whose Sasaki-Reeb cone increases in dimension. The following remark exhibits the procedure to obtain the aforementioned perturbations.
\medskip


\begin{rmk}
Let us consider the orbifold  $X_{f}\subset \mathbb P(w_0, w_1, w_2, w_3, w_4)$  cut out by the degree $d$ polynomial of chain-chain type  of the form $$f=z_{0}^{a_{0}}+z_{0}z_{1}^{2}+z_{2}^{a_{2}}+z_{2}z_{3}^{a_{3}}+z_{3}z_{4}^{2}.$$ 
Then there is a perturbation  of  the Berglund-Hübsch dual 
 $f^T=z_0^{a_0} z_1+z_1^2+z_2^{a_2} z_3+z_3^{a_3} z_4+z_4^2,$ a  polynomial of the form 
$$g=z_{0}^{\tilde{a}_0}+z_{1}^{2}+z_{2}^{a_{2}}z_{3}+z_{3}^{\tilde{a}_{3}}+z_{4}^{2}\in H^0(\mathbb{P}(\tilde{\mathbf{w}}), \mathcal{O}(\tilde{d})),$$  
with  $\tilde{\mathbf{w}}=(\tilde{w}_0, \tilde{w}_1, \tilde{w}_2, \tilde{w}_3, \tilde{w}_4)$ and degree $\tilde{d}$  determined by the  Berglund-H\"ubsch transpose $f^T$ of $f,$ such that  $g$ determines a link whose Sasaki-Reeb cone has dimension 2.
Indeed, the matrix of exponents $A$ of $f$ given by 
\begin{equation*}
    A =\begin{bmatrix}
      a_0 & 0 & 0 & 0 & 0 \\
      1 & 2 & 0 & 0 & 0 \\
      0 & 0 & a_2 & 0 & 0 \\
      0 & 0 & 1 & a_3 & 0 \\
      0 &  0 & 0 & 1 & 2
    \end{bmatrix}
\end{equation*}
satisfies the matrix equation $$AW=D$$ with $W=(w_0,w_1, w_2, w_3,w_4)^T$ and 
$D=(d,d,d,d,d)^T.$ From this we obtain: 
\begin{equation}
    a_0w_0=d,\quad a_2w_2=d,
\end{equation}
\begin{equation}
    w_2+a_3w_3=d.
\end{equation}

For the Berglund-Hübsch transpose $f^T$ we have the following matrix equation 

\begin{equation*}
    \begin{bmatrix}
      a_0 & 1 & 0 & 0 & 0 \\
      0 & 2 & 0 & 0 & 0 \\
      0 & 0 & a_2 & 1 & 0 \\
      0 & 0 & 0 & a_3 & 1 \\
      0 &  0 & 0 & 0 & 2
    \end{bmatrix}
    \begin{bmatrix}
    \tilde{w}_0\\
    \tilde{w}_1\\
    \tilde{w}_2\\
    \tilde{w}_3\\
    \tilde{w}_4\\
    \end{bmatrix}
    =\begin{bmatrix}
       \tilde{d}\\
       \tilde{d}\\
       \tilde{d}\\
       \tilde{d}\\
       \tilde{d}\\
    \end{bmatrix}.
\end{equation*}
A solution to this system is given by the following equations

\begin{equation}
    \tilde{w}_0=w_0(d-w_2),
\end{equation}
\begin{equation}
    \tilde{w}_1=d(d-w_2),
\end{equation}
\begin{equation}
    \tilde{w}_2=w_2(2d-2w_2-w_3),
\end{equation}
\begin{equation}
    \tilde{w}_3=dw_3,
\end{equation}
\begin{equation}
    \tilde{w}_4=d(d-w_2),
\end{equation}
\begin{equation}
    \tilde{d}=2d(d-w_2).
\end{equation}
By the first equality given in equation (5) and  equation (6), we have $w_{0}\mid d$ and $w_{3}\mid (d-w_{2})$. From equations (5), (8), (11) and (12) we obtain the polynomial 
$$g=z_{0}^{\tilde{a}_0}+z_{1}^{2}+z_{2}^{a_{2}}z_{3}+z_{3}^{\tilde{a}_{3}}+z_{4}^{2}\in H^0(\mathbb{P}(\tilde{\mathbf{w}}), \mathcal{O}(\tilde{d})),$$  
where $\tilde{a}_i=2a_i$ for $i=1,3.$ It follows that  the connected component of the Sasaki automorphism group for the polynomial $g$ is 
$U(1) \times S O(k)$ where $k$ is the number of uncoupled quadratic terms found in the polynomial which in this case is either $2$ or $3$. Thus,  the Sasaki-Reeb cone has dimension 
$\left\lfloor \frac{k}{2}\right\rfloor +1=2$ in this case. (See also Theorem 4.1 below.) 
\end{rmk}


\section{Extremal Sasaki metrics and Berglund-Hübsch transpose}
In this section we combine a precise and explicit notion of K-stability developed by Boyer and van Coevering in \cite{BvC} and the Berglund-Hübsch transpose rule to produce examples  motivated by the questions given in the introduction: 
\begin{itemize}
\item Are there invertible polynomials besides Brieskorn-Pham polynomials whose corresponding links are homotopy spheres and admit extremal metrics?
\item  Under what conditions is the emptiness of the extremal Sasaki-Reeb cone  stable via the Berglund-H\"ubsch transpose rule?  
\end{itemize}

In Subsection 4.1 we produce examples that provides some mild evidence that suggests that the answer to the first question may be no.  In \cite{BvC}, Boyer and van Coevering defined a version of the K-stability of Collins and Székelyhidi \cite{CS} relative to a maximal torus and showed that the Lichnerowicz obstruction \cite{GMSY}, which obstructs the existence of Sasaki-Einstein metrics for the Reeb vector field associated to the embedding, can be improved to  obtain a {\it generalized Lichnerowicz obstruction}: that is,  an obstruction to  the existence of extremal Sasaki metrics for every Reeb vector field in the Sasaki-Reeb cone, see Theorem 4.1 below. This result  allows us to produce families of Fano cone singularities $X$ of complexity 3 or complexity 4 such that their corresponding links $L=X \cap S^{2 N-1}$ of positive Sasakian metric do not admit any Sasaki extremal metrics  in the whole Sasaki-Reeb cone. We manufacture these examples from Berglund-H\"ubsch duals  of links of  invertible polynomials with chain-cycle type singularities of the form $$f=z_{0}^{a_{0}}+z_{0}z_{1}^{a_{1}}+z_{4}z_{2}^{a_{2}}+z_{2}z_{3}^{a_{3}}+z_{3}z_{4}^{a_{4}}$$ which were studied in \cite{CGL}. These  duals are  not well-formed links and under certain conditions on the weights and the degree they admit Sasaki-Einstein metrics. From these links, we obtain homotopy 9-spheres with Sasaki-Reeb cone of dimension 2 and these are given as 2-fold branched covers of $S^9$  whose branching loci are rational homology 7-spheres arising as  links of  not well-formed orbifolds described by invertible polynomials described as Thom-Sebastiani sum of singularities of  chain and cycle type, see Theorem 4.3. Similarly, we produce examples of rational homology $(4k+3)$-spheres and  homotopy $(4k+1)$-spheres with Sasaki-Reeb cone of dimension $k$  with no extremal Sasaki metrics. All these examples spheres produced here are links of invertible 
polynomials that cannot be written as Brieskorn-Phams polynomials.

With respect to the second question, in Subsections 4.2 and 4.3, we exhibit families of links $L_f(\mathbf{w}, d)$ with one dimensional Sasaki-Reeb cone not admitting extremal Sasaki metrics that via the  Berglund-H\"ubsch rule produce  links with Sasaki-Reeb cones of dimension greater than one and such that their whole Sasaki-Reeb cones do not admit extremal Sasaki metrics. More precisely, the dual link $L_{f^{T}}(\tilde{\mathbf{w}}, \tilde{d})$ admits a perturbation in its local moduli, that is,  an element in  the complex space of weighted homogeneous polynomials of degree $\tilde{d}$,  $H^0(\mathbb{P}(\tilde{\mathbf{w}}),\mathcal{O}(\tilde{d})),$ which determines a link whose Sasaki-Reeb cone is enlarged and  which remains obstructed to admitting  extremal Sasaki metrics.  

From Theorems 4.4, 4.6 and  4.8, and Lemmas 4.9-4.12, we notice that under some arithmetic constraints,  
the Berglund-Hübsch transpose rule can be interpreted as a map that transform the classical Lichnerowicz obstruction into the improved version of this obstruction for Thom-Sebastiani sums involving chain type singularities. It would be interesting to determine, in a systematic manner, a larger family of invertible polynomials with that feature. 

We recall the Lichnerowicz obstruction: 
\medskip

\noindent{\bf Lichnerowicz obstruction (\cite{GMSY}).}  
{\it Let $L_{f}$ be a smooth link of an isolated hypersurface singularity associated to a polynomial $f$ of $n+1$ variables with weight vector $\bf w $ and  degree $d$. If the index $I$ satisfies the inequality 
    \begin{equation*}
        I=|\mathbf{w}|-d > n\min_{i}w_{i},
    \end{equation*}
    then the link $L_{f}$ does not admit a Sasaki-Einstein metric.}
\medskip

Now, let us consider a weighted homogeneous polynomial  $f:\mathbb{C}^{n+1}\rightarrow \mathbb{C}$ of degree $d$ and weight vector $\mathbf{w}=\left(w_0, \ldots, w_n\right)$ with only an isolated singularity at the origin such that  $f$ is given of the form
$$
f\left(z_0, \ldots, z_n\right)=f^{\prime}\left(z_0, \ldots, z_k\right)+z_{k+1}^2+\cdots+z_n^2
$$
with $n-k \geq 2$ and all weights $w_i$ with $i=0, \ldots, k$ satisfying $2 w_i<d^{\prime}$, the degree of $f^{\prime}$. 
In this case the connected component of the Sasaki automorphism group is $U(1) \times S O(n-k)$ so the Sasaki-Reeb cone has dimension $\left\lfloor \frac{n-k}{2}\right\rfloor +1.$ 
Let us assume no linear factor is allowed otherwise the link would be a standard sphere. For this sort of polynomials one can extend the obstruction  given above to a generalized version of the Lichnerowicz obstruction for existence of extremal Sasaki metrics:

\begin{thm}[\cite{BvC}]
Under the hypotheses given above,  the dimension of the  Sasaki-Reeb cone of the link $L_f$ is given by $\left\lfloor \frac{n-k}{2}\right\rfloor +1.$ 
Furthermore, if the generalized Lichnerowicz obstruction holds, that is, if the following inequality holds 
\begin{equation}
\sum_{i=0}^k w_i-\min_l\{w_l\}n +\frac{d}{2}(n-k-2) \geq 0, 
\end{equation} 
then there are no extremal Sasaki metrics in the entire Sasaki-Reeb cone of the link $L_f$.
\end{thm}

Boyer and van Coevering applied these  results to  Brieskorn-Pham polynomials to exhibit the first examples of Sasaki manifolds with Sasaki-Reeb cones of dimension greater than one which do not admit extremal metrics.  Our approach illustrates how one can get a lot more examples if one uses invertible polynomials together with the Berglund-H\"ubsch transpose rule.

From this powerful criterion, we have the following observation with respect to Remark 3.1. 

\begin{lemma}
Let us consider the orbifold  $X_{f}\subset \mathbb P(w_0, w_1, w_2, w_3, w_4)$  with $w_0+w_3\geq 6w_2,$ cut out by the degree $d$ chain-chain type  polynomial of the form $$f=z_{0}^{a_{0}}+z_{0}z_{1}^{2}+z_{2}^{a_{2}}+z_{2}z_{3}^{a_{3}}+z_{3}z_{4}^{2}.$$ 
Then the polynomial  
$$g=z_{0}^{\tilde{a}_0}+z_{1}^{2}+z_{2}^{a_{2}}z_{3}+z_{3}^{\tilde{a}_{3}}+z_{4}^{2}\in H^0(\mathbb{P}(\tilde{\mathbf{w}}), \mathcal{O}(\tilde{d})),$$ which is the perturbation of the Berglund-H\"ubsch dual $f^T$  described in Remark 3.1, is such that  there are no extremal Sasaki metrics in the whole  Sasaki-Reeb cone of dimension $2$ for the associated link $L_{g}$.
\end{lemma}
\begin{pf}

We have to show the inequality (13), that is, we need to verify  
$$\tilde{w}_0+\tilde{w}_2+\tilde{w}_3-4\min_{l}\{\tilde{w}_l\}\geq 0.$$ The only  candidates for minimum are $\tilde{w}_0,$ $\tilde{w}_2,$ and $\tilde{w}_3.$ We have 
\begin{eqnarray*}
  \tilde{w}_0+\tilde{w}_2+\tilde{w}_3 &=& w_0(d-w_2)+w_2(2d-2w_2-w_3)+dw_3\\
  &=& (d-w_2)(w_0+2w_2+w_3)\\
  &\geq & 8(d-w_2)w_2.
\end{eqnarray*}
The last inequality follows from the hypothesis on the weights. 
Notice $$2(d-w_2)w_2\geq w_2(2d-2w_2-w_3),$$ thus $8(d-w_2)w_2\geq 4 \tilde{w}_2$ which proves the lemma.
\hfill$\square$
\end{pf}
\medskip

\subsection{Examples derived from chain-cycle type singularities} In this subsection, 
we use  Theorem 4.1 on certain invertible polynomials that can be written as a Thom-Sebastiani sum of a  Brieskorn-Pham polynomial and a cycle type polynomial and are obtained under the Berglund-H\"ubsch transpose from chain-cycle type singularities. Actually from the list of rational homology 7-spheres admitting Sasaki-Einstein metrics given in \cite{CGL} we produce  37 new examples in the form of 31 homotopy 9-spheres and 6 examples of manifolds of the form $S^4\times S^5$, such that its Sasaki-Reeb cones of dimension 2 do not admit extremal Sasaki metrics. 
We can extrapolate the methods to exhibit examples of rational homology 
$(4k+3)$-spheres and  homotopy $(4k+1)$-spheres with Sasaki-Reeb cone of dimension  $k$  with no extremal Sasaki metrics. In particular, we produce Fano cone singularities of complexity 3 or complexity 4.

First we collect some facts on certain links of chain-cycle singularities. This will be helpful to find the Alexander polynomial of $f^T.$ Consider the  hypersurface $X_{f}$ where $f$ is a degree $d$ polynomial of chain-cycle type of the form 
$$f=z_{0}^{a_{0}}+z_{0}z_{1}^{a_{1}}+z_{4}z_{2}^{a_{2}}+z_{2}z_{3}^{a_{3}}+z_{3}z_{4}^{a_{4}}$$ with $\mathbf{w}=(w_{0},w_{1},w_{2},w_{3},w_{4})=(m_{3}v_{0},m_{3}v_{1},m_{2}v_{2},m_{2}v_{3},m_{2}v_{4}),$ such that  $\gcd(m_3,m_2)=1$ and $d=m_3m_2.$ It is not difficult to obtain for the chain-cycle polynomial given above that 
$$v_{0}=\frac{w_0}{\gcd(w_0, d)}=1$$ which implies that $w_0=m_3.$ Thus, we obtain
\begin{equation*}
            w_{0}=\dfrac{d}{a_{0}}, w_{1}=\dfrac{d-w_{0}}{a_{1}}
    \end{equation*} 
    \begin{equation*}
        a_{0}=\dfrac{d}{w_{0}}=m_{2} \mbox{ and } a_{1}=\dfrac{d-w_{0}}{w_{1}}=\dfrac{m_{3}m_{2}-m_{3}}{m_{3}v_{1}}=\dfrac{m_{2}-1}{v_{1}}.
    \end{equation*}

Using the Berglund-Hübsch transpose rule, we have its corresponding  polynomial $${f}^{T}=z_{0}^{a_{0}}z_{1}+z_{1}^{a_{1}}+z_{3}z_{2}^{a_{2}}+z_{4}z_{3}^{a_{3}}+z_{2}z_{4}^{a_{4}}$$ with weights $\tilde{\mathbf{w}}=(\Tilde{w}_{0},\Tilde{w}_{1},\Tilde{w}_{2},\Tilde{w}_{3},\Tilde{w}_{4})$ and degree $\Tilde{d}.$ 

Consider  $\Tilde{u}_{i}=\dfrac{\Tilde{d}}{\gcd(\Tilde{d},\Tilde{w}_{i})}\quad \mbox{and}\quad  \Tilde{v}_{i}=\dfrac{\Tilde{w}_{i}}{\gcd(\Tilde{d},\Tilde{w}_{i})}.$ 
We have the following facts, whose proofs can be found in the proof of Theorem 3.2 in  \cite{CGL}:
\begin{enumerate}
    \item[{\bf Fact 1}] The degree of $f^T$ equals $$\Tilde{d}=d(m_{2}-1)=m_{3}m_{2}(m_{2}-1)$$ and the index $\tilde{I}$ is given by  
    $$\tilde{I}=\vert \tilde{\mathbf{w}}\vert -\tilde{d}=m_{2}-1.$$ 
    \item[{\bf Fact 2}] 
Furthermore, we have    $\Tilde{u}_{0}=m_{2}a_{1}, \ \Tilde{u}_{1}=a_{1}, \ \Tilde{u}_{2}=\Tilde{u}_{3}=\Tilde{u}_{4}=m_{3}$
and $\Tilde{v}_{0}=a_{1}-1, \  \Tilde{v}_{1}=1, \ \Tilde{v}_{2}=1-a_{4}+a_{3}a_{4}, \ \Tilde{v}_{3}=1-a_{2}+a_{4}a_{2}, \ \Tilde{v}_{4}=1-a_{3}+a_{2}a_{3}.$
    \item[{\bf Fact 3}]  The third Betti number is given by $$b_3(L_{f^T})=-1+\frac{a_2a_3a_4+1}{m_3}.$$ 
    \item[{\bf Fact 4}] In particular, for $L_{f^T}$ a rational homology 7-sphere, we have $m_3=a_2a_3a_4+1.$  
\begin{enumerate}
\item From equality (4) we have
{\small{
\begin{equation}
\left(\dfrac{1}{\Tilde{v}_{2}}+\dfrac{1}{\Tilde{v}_{3}}+\dfrac{1}{\Tilde{v}_{4}}\right)-m_{3}\left(\dfrac{1}{\Tilde{v}_{2}\Tilde{v}_{3}}+\dfrac{1}{\Tilde{v}_{3}\Tilde{v}_{4}}+\dfrac{1}{\Tilde{v}_{4}\Tilde{v}_{2}}\right)+\dfrac{\Tilde{m}_{3}^{2}}{\Tilde{v}_{2}\Tilde{v}_{3}\Tilde{v}_{4}} =1.
\end{equation}
}}

\item The weights are given by 
$$\tilde{w_{1}}=\dfrac{\tilde{d}}{a_{1}}=\dfrac{m_{3}m_{2}(m_{2}-1)}{a_{1}}$$
$$\tilde{w_{0}}=\dfrac{\tilde{d}-\tilde{w_{1}}}{a_{0}}=\dfrac{m_{3}m_{2}(m_{2}-1)-m_{3}m_{2}v_1}{m_{2}}=\dfrac{m_{3}(m_{2}-1)(a_{1}-1)}{a_{1}}$$
$$\tilde{w}_{2}=\dfrac{\tilde{d}(1-a_{4}+a_{4}a_{3})}{1+a_{2}a_{3}a_{4}}=m_{2}(m_{2}-1)(1-a_{4}+a_{3}a_{4}) = m_{2}(m_{2}-1)\tilde{v}_{2}$$
$$\tilde{w}_{3}=\dfrac{\tilde{d}(1-a_{2}+a_{2}a_{4})}{1+a_{2}a_{3}a_{4}}=m_{2}(m_{2}-1)(1-a_{2}+a_{2}a_{4})= m_{2}(m_{2}-1)\tilde{v}_{3}$$
$$\tilde{w}_{4}=\dfrac{\tilde{d}(1-a_{3}+a_{2}a_{3})}{1+a_{2}a_{3}a_{4}}=m_{2}(m_{2}-1)(1-a_{3}+a_{2}a_{3}) = m_{2}(m_{2}-1)\tilde{v}_{4}.$$
In particular, for $a_{1}=2$, then the formulas for $\tilde{w}_{0}$ and $\tilde{w}_{1}$ can be written as:

$$\tilde{w}_{0}=\dfrac{m_{3}(m_{2}-1)}{2}  \ \ \ 
 \mbox{ and } \ \ \ \ 
 \tilde{w}_{1}=\dfrac{m_{3}m_{2}(m_{2}-1)}{2}.$$

\item There are three possible values for $H_{3}(L_{f^T},\mathbb{Z})$:
\begin{itemize}
   \item   If $\gcd(a_{1},m_{3})-1<1$, then $H_{3}(L_{f^T},\mathbb{Z})=\mathbb{Z}_{m_{3}}$.
   \item  If $\gcd(a_{1},m_{3})-1=1$, then $H_{3}(L_{f^T},\mathbb{Z})=\mathbb{Z}_{d}$.
   \item  If $\gcd(a_{1},m_{3})-1>1$, then $H_{3}(L_{f^T},\mathbb{Z})=\mathbb{Z}_{d}\oplus \underbrace{\mathbb{Z}_{m_{2}}\oplus\cdots\oplus\mathbb{Z}_{m_{2}}}_{(\gcd(a_{1},m_{3})-2)-times}$.
 \end{itemize}
    
\item Additionally, if the polynomial $g$ is of  the form 
\begin{eqnarray*}
g &=&  f^{T}+z_{5}^{2}+z_{6}^{2}+\ldots +z_{n}^{2}\\
&=& z_{0}^{a_{0}}z_{1}+z_{1}^{2}+z_{3}z_{2}^{a_{2}}+z_{4}z_{3}^{a_{3}}+z_{2}z_{4}^{a_{4}}+z_{5}^{2}+z_{6}^{2}+\ldots +z_{n}^{2},
\end{eqnarray*}
One can find a perturbation of $g$, a  polynomial  of the form $$g_2=z_{0}^{\tilde{a}_{0}}+z_{1}^{2}+z_{3}z_{2}^{a_{2}}+z_{4}z_{3}^{a_{3}}+z_{2}z_{4}^{a_{4}}+z_{5}^{2}+z_{6}^{2}+\ldots +z_{n}^{2}$$ 
with weight vector $\textbf w_{Ext}=(\tilde{\mathbf{w}},w_{5},w_{6}, \ldots , w_{n})$ and degree $\tilde{d}$ 
where $\tilde{a}_{0}=2m_2$. This follows from a similar argument as the one given in Remark 3.1 and 
the fact that $\tilde{d}=m_{3}m_{2}(m_{2}-1)$ is divisible by $\tilde{w}_{0}=\dfrac{m_{3}(m_{2}-1)}{2}.$

\end{enumerate}

\end{enumerate}

Notice that under certain circumstances, we can simplify the values for $\tilde{w}_{i}'s$ and degree $\tilde{d}$. Nevertheless, since we know these values explicitly, we would rather work with the formulas obtained above. In the next theorem, the conditions on the degrees and weight vectors of the corresponding polynomials forbid them to be expressed as  Brieskorn-Pham polynomials.

\begin{thm}
     Consider a polynomial of chain-cycle type $$f=z_{0}^{a_{0}}+z_{0}z_{1}^{a_{1}}+z_{4}z_{2}^{a_{2}}+z_{2}z_{3}^{a_{3}}+z_{3}z_{4}^{a_{4}}$$ with $a_1=2$  that  cuts out a projective hypersurface of degree $d$ in weighted projective space $\mathbb{P}(w_0, w_1, \ldots , w_4)$ such that $$(w_{0},w_{1},w_{2},w_{3},w_{4})=(m_{3}v_{0},m_{3}v_{1},m_{2}v_{2},m_{2}v_{3},m_{2}v_{4}),$$ with $m_3$ odd, $\gcd(m_{2},m_{3})=1$ and $d=m_3m_2.$ Also, suppose that $f$ determines a link $L_f$ with the rational homology sphere.
Then the polynomial $$g=f^T+z_5^2+\ldots +z_n^2$$ with $f^T$ the Berglund-Hübsch dual of $f,$  
determines a link $L_{g_2}$ with $g_2$ given as in Item 4(d) above: $$g_2=z_{0}^{\tilde{a}_{0}}+z_{1}^{2}+z_{3}z_{2}^{a_{2}}+z_{4}z_{3}^{a_{3}}+z_{2}z_{4}^{a_{4}}+z_{5}^{2}+z_{6}^{2}+\ldots +z_{n}^{2},$$ 
where $\tilde{a}_{0}=2m_2$,  such that 
\begin{enumerate}
     \item[(1)] The Sasaki-Reeb cone of $L_{g_2}$ has dimension $1+\lfloor{\frac{n-3}{2}}\rfloor$ and there are no extremal Sasaki metrics in its whole Sasaki-Reeb cone. 
    \item[(2)]  If $n$ is even, then $L_{g_2}$ is a rational homology $(2n-1)$-sphere. 
    \item[(3)]  If $n$ is odd, then  $L_{g_2}$ is a  homotopy $(2n-1)$-sphere and $\Delta_g(-1)=m_3.$ In particular the diffeomorphism type of $L_{g_2}$ is determined by $m_3.$
    
\end{enumerate}  
\end{thm}
\begin{pf} 
First, notice that, since $\tilde{w}_{0}=\dfrac{m_{3}(m_{2}-1)}{2}$,  $m_{2}$ must be odd. The first statement of Item (1) follows from Theorem 4.1. For the second statement of Item (1),  we need to prove that
\begin{equation}
\tilde{w}_{0}+\tilde{w}_{2}+\tilde{w}_{3}+\tilde{w}_{4}-n\min_{l}\{\tilde{w}_{l}\}+(n-5)\dfrac{\tilde{d}}{2}\geq0 \,\, \mathrm{for\,\, all}\, n\geq 6.
\end{equation}
We use induction on $n$. For $n=6$, the inequality above is
$$\tilde{w}_{0}+\tilde{w}_{2}+\tilde{w}_{3}+\tilde{w}_{4}-6\min_{l}\{\tilde{w}_{l}\}+\dfrac{\tilde{d}}{2}\geq0$$
Since $\tilde{w}_{0}+\tilde{w}_{2}+\tilde{w}_{3}+\tilde{w}_{4}\geq4\min_{l}\{\tilde{w}_{l}\}$, it is sufficient to show that $\dfrac{\tilde{d}}{2}\geq 2\min_{l}\{\tilde{w}_{l}\}$. As $m_{2}\geq 2$, then
$$\tilde{d}=m_{3}m_{2}(m_{2}-1)\geq 2m_{3}(m_{2}-1)= 4\left( \dfrac{m_{3}(m_{2}-1)}{2}\right) =4\tilde{w}_{0}.$$
Therefore, we have
\begin{equation}
\dfrac{\tilde{d}}{2}\geq 2\tilde{w}_{0}\geq 2\min_{l}\{\tilde{w}_{l}\}.
\end{equation}
Thus, the inequality is true for $n=6$.

On the other hand, let us suppose that  the inequality (15) holds for some $n\geq6$. Now, we are going to show this for $n+1$. As $\tilde{d}=2w_n$ we have  
\begin{equation}
\dfrac{\tilde{d}}{2}- \min_{l}\{\tilde{w}_{l}\}\geq 0.
\end{equation}
Adding (15) and  (17) we obtain

\begin{equation*}
    \tilde{w}_{0}+\tilde{w}_{2}+\tilde{w}_{3}+\tilde{w}_{4}-(n+1)\min_{l}\{\tilde{w}_{l}\}+((n+1)-5)\dfrac{\tilde{d}}{2} \geq 0.
\end{equation*}

Next, we will prove the result above for $n=5$. 
\begin{enumerate}
\item[a)] Let $m_{2}\geq 5$.
    In this case, we will show that
    $$w_{0}+w_{2}+w_{3}+w_{4}-5\min_{l}\{w_{l}\}\geq 0.$$
First, we suppose that $w_{0}=\min_{l}\{w_{l}\}$, then we have
   $$ w_{0}+w_{2}+w_{3}+w_{4}-5\min_{l}\{w_{l}\} = \tilde{d}+\tilde{I}-\tilde{w}_{1}-5\tilde{w}_{0},$$ that  can be written as  
    $$m_{3}m_{2}(m_{2}-1) + (m_{2}-1)-\dfrac{m_{3}m_{2}(m_{2}-1)}{2}-\dfrac{5m_{3}(m_{2}-1)}{2}$$    
    which equals
    $$( m_{2}-1)\left( \dfrac{m_{3}}{2}(m_{2}-5)+1\right).$$
Since $m_{2}\geq 5$, we obtain $w_{0}+w_{2}+w_{3}+w_{4}-5\min_{l}\{w_{l}\}>0$.

On the other hand, if $\tilde{w_{i}}=\min_{l}\{w_{l}\}$ for $i=2,3$ or $4$, we can suppose without loss of generality that $\tilde{w_{2}}=\min_{l}\{w_{l}\}$. Therefore, we have
$\tilde{w}_{0}>\tilde{w}_{2}$ implies  $m_{3}>2m_{2}\tilde{v}_{2}.$
Then
\begin{align*}
    w_{0}+w_{2}+w_{3}+w_{4}-5\min_{l}\{w_{l}\} & = \tilde{d}+\tilde{I}-\tilde{w}_{1}-5\tilde{w}_{2} \\
    & = ( m_{2}-1)\left( \dfrac{m_{3}m_{2}}{2}+1 -5m_{2}\tilde{v}_{2}\right) \\
    & > ( m_{2}-1)\left( m_{2}^{2}\tilde{v}_{2}+1 -5m_{2}\tilde{v}_{2}\right) \\
    & = ( m_{2}-1)\left( m_{2}\tilde{v}_{2}(m_{2}-5)+1\right).
\end{align*}

As $m_{2}\geq5$, we obtain $$w_{0}+w_{2}+w_{3}+w_{4}-5\min_{l}\{w_{l}\}>0.$$
\item[b)] Let $m_{2}=3.$ In this case, the weights are given by
$$\tilde{w}_{0}=m_3, \ \ \tilde{w}_{1}=3m_{3}, \ \ \tilde{w}_{2}=6\tilde{v}_{2}, \ \ \tilde{w}_{3}=6\tilde{v}_{3}, \ \ \mbox{ and } \ \ \tilde{w}_{4}=6\tilde{v}_{4}.$$
In addition, the degree and the index are $\tilde{d}=6m_{3}$ and $\tilde{I}=2$, respectively. Since the weights $\tilde{w}_{2}, \tilde{w}_{3}$ and $\tilde{w}_{4}$ come from cycle block of polynomial $g$, 
we can suppose, without loss of generality, that $\tilde{w}_{4}> \tilde{w}_{0}$.

On the other hand, since $\sum_{i=0}^{4}\tilde{w}_{i}=\tilde{d}+\tilde{I}$, we have 
$$\tilde{w}_{2}+\tilde{w}_{3}+\tilde{w}_{4}=2m_{3}+2.$$
Since $\tilde{w}_{4}>\tilde{w}_{0}=m_{3}$, we have $\tilde{w}_{2}+\tilde{w}_{3}<m_{3}+2$. It implies $\tilde{w}_{2}\leq \frac{m_{3}+1}{2}$ or $\tilde{w}_{3}\leq \frac{m_{3}+1}{2}$. Thus, we obtain

\begin{align*}
    w_{0}+w_{2}+w_{3}+w_{4}-5\min_{l}\{w_{l}\} & = \tilde{d}+\tilde{I}-\tilde{w}_{1}-5\min_{l}\{w_{l}\} \\
    & = 3m_{3}+2-5\min_{l}\{w_{l}\}\\
    & \geq 3m_{3}+2-5\left(\dfrac{m_{3}+1}{2}\right) \\
    & = \dfrac{m_{3}-1}{2}.
\end{align*}
As $m_{3}\geq1$, we obtain $w_{0}+w_{2}+w_{3}+w_{4}-5\min_{l}\{w_{l}\}\geq 0.$
\end{enumerate}

For  Item (2), the case $n$ even follows from Corollary 9.5.3 in \cite{BBG} and Theorem 4.1 in \cite{CGL}. For Item (3), notice this  reduces to the polynomial $g=f^T+z_5^2$, since $(\Lambda_2-\Lambda_1)(\Lambda_2-\Lambda_1)=\Lambda_1.$ From Equation (1) we have that the Alexander polynomial for $g$ is given by 
\begin{equation}
\operatorname{div}\Delta_g=\prod_{i=0}^5\left (  \frac{\Lambda_{\tilde{u}_i}}{\tilde{v}_i} -\Lambda_1 \right ),
\end{equation}
 where the ${\tilde{u}_i}$'s and the ${\tilde{v}_i}$'s are given as before. In particular for the given  data we have 
${\tilde{u}_0}=2m_2,$ ${\tilde{u}_1}={\tilde{u}_5}=2$ and ${\tilde{u}_2}={\tilde{u}_3}={\tilde{u}_4}=m_3$
and ${\tilde{v}_0}={\tilde{v}_1}={\tilde{v}_2}=1.$ Hence, from Equation (1), using the relation $\Lambda_a \Lambda_b=\operatorname{gcd}(a, b) \Lambda_{\operatorname{lcm}(a, b)}$ several times,  we have 
{\footnotesize{
\begin{eqnarray*}
   \operatorname{div}\Delta_g &=& \left (\Lambda_{2m_2} -\Lambda_1 \right )
   \left (\Lambda_2 -\Lambda_1 \right )\left (\frac{\Lambda_{\tilde{u}_2}}{\tilde{v}_2} -\Lambda_1 \right )\left (\frac{\Lambda_{\tilde{u}_3}}{\tilde{v}_3} -\Lambda_1 \right )\left (\frac{\Lambda_{\tilde{u}_4}}{\tilde{v}_4} -\Lambda_1 \right )\left (\Lambda_2 -\Lambda_1 \right )\\
   &=& \left (\Lambda_{2m_2} -\Lambda_1 \right )\left (\frac{\Lambda_{m_3}}{\tilde{v}_2} -\Lambda_1 \right )\left (\frac{\Lambda_{m_3}}{\tilde{v}_3} -\Lambda_1 \right )\left (\frac{\Lambda_{m_3}}{\tilde{v}_4} -\Lambda_1 \right )\\
   &=& \left (\Lambda_{2m_2} -\Lambda_1 \right )\left (\Lambda_{m_3}
   \left \{ \frac{m_3^2}{\tilde{v}_2\tilde{v}_3\tilde{v}_4} 
   -m_3\left (\frac{1}{\tilde{v}_2\tilde{v}_3}+\frac{1}{\tilde{v}_4\tilde{v}_3}+\frac{1}{\tilde{v}_2\tilde{v}_3}
   \right) + \left ( \frac{1}{\tilde{v}_2}+\frac{1}{\tilde{v}_3}+\frac{1}{\tilde{v}_4}  \right )  \right \}-\Lambda_1 \right ).\\
\end{eqnarray*}
}}
From Theorem 3.2 in \cite{CGL} it follows that $L_{f^T}$ is a rational homology sphere, then using Equation (4)  it follows that  
\begin{eqnarray*}
\operatorname{div}\Delta_g &=& \left (\Lambda_{2m_2} -\Lambda_1 \right )\left (\Lambda_{m_3} -\Lambda_1 \right )\\
&=& \Lambda_{2m_2}\Lambda_{m_3}-\Lambda_{m_3}-\Lambda_{2m_2}+\Lambda_1\\
&=& \Lambda_{2d}-\Lambda_{m_3}-\Lambda_{2m_2}+1.
\end{eqnarray*}
The third equality above  follows from $\gcd(2m_2, m_3)=1$ since $m_3$ is odd. Thus the Alexander polynomial is given by 
\begin{eqnarray*}
\Delta_g(t)=\frac{(t^{2d}-1)(t-1)}{(t^{m_3}-1)(t^{2m_2}-1)}.
\end{eqnarray*}
It follows that $\Delta_g(1)=\frac{2d}{m_3(2m_2)}=1.$ Thus $L_g$ is a homotopy 9-sphere. Furthermore, performing the change of variable $t^{2m_2}=l,$ we have 
\begin{eqnarray*}
\Delta_g &=&\frac{(l^{m_3}-1)(t-1)}{(t^{m_3}-1)(l-1)}\\
&=& \frac{(l^{m_3-1}+\ldots +1)}{(t^{m_3-1}\ldots +1)}.
\end{eqnarray*}
Hence $\Delta_g(-1)=m_3.$
\hfill$\square$
\end{pf}

 \begin{rmk}
With respect to Item (2) in Theorem 4.3, we can say a bit more. We explain this for the simplest case $h=f^{T}+z_{5}^{2}+z_{6}^{2},$  we claim that  $H_{3}(L_{f},\mathbb{Z})$ and $H_{5}(L_{h},\mathbb{Z})$ have the same torsion. 
    For this, we use the Orlik's algorithm. Computing the values $u_{i}'s$ and $v_{i}'s$:
\begin{center}
      \begin{tabular}{| c || c | c | c | c | c | c | c |}
\hline
$i$ & $0$ & $1$ & $2$ & $3$ & $4$ & $5$ & $6$\\ \hline
$u_{i}$ &  $2m_{2}$ & $2$ & $m_{3}$ & $m_{3}$ & $m_{3}$ & $2$ & $2$\\ \hline
$v_{i}$ & $1$ & $1$ & $\tilde{v}_{2}$ & $\tilde{v}_{3}$ & $\tilde{v}_{4}$ & $1$ & $1$
\\ \hline
\end{tabular}
\end{center}
Next, we calculate the numbers $c_{i_{1}\dots i_{s}}$. For the polynomial $h$, we obtain $c_{\emptyset}=\gcd(2,m_{3})$, $c_{234}=\dfrac{2}{\gcd(2,m_{3})}$, $c_{0156}=\dfrac{m_{3}}{\gcd(2,m_{3})}, c_{123456}=m_{2}$ and $c_{i_{1}\dots i_{s}}=1$ in other cases. Thus, we only need to calculate the values $k_{i_{1}\dots i_{s}}$ where $c_{i_{1}\dots i_{s}}\neq 1$. Therefore, we obtain:
$$k_{\emptyset}=1, \ \ \ \ k_{234}=0, \ \ \ \  k_{0156}=1, \ \ \  \ \mbox{and} \ \ \ \ k_{123456}=\gcd(2,m_{3})-1$$

Then, there are two cases for $H_{5}(L_{h},\mathbb{Z})$:
\begin{itemize}
    \item If $\gcd(2,m_{3})-1<1$, then $H_{5}(L_{h},\mathbb{Z})=\mathbb{Z}_{m_{3}}$.
    \item If $\gcd(2,m_{3})-1=1$, then $H_{5}(L_{h},\mathbb{Z})=\mathbb{Z}_{m_{3}m_{2}}$.
\end{itemize}
It follows from Fact 4(c) that the torsion of $H_{5}(L_{h},\mathbb{Z})$ coincides with the torsion of $H_{3}(L_{f},\mathbb{Z})$ when $a_{1}=2$.
\end{rmk}

Now, we present a representative example of Theorem 4.3.
\begin{example}
Consider the rational homology 7-sphere $$L_f(\mathbf{w},d)=f^{-1}(0)\cap S^9$$ with weight vector  $\mathbf{w}=(701,701,198,381, 123)$ and degree $d=2103.$ It was proven in \cite{CL} that this link admits Sasaki-Einstein metric. 
One can describe $f$ as an invertible polynomial of BP-cycle type, cycle-cycle type and 
chain-cycle type. We are interested in a representation given by the latter type where the Berglund-Hübsch transpose rule acts in a non-trivial fashion, as shown in \cite{CGL}. 
Explicitly, let us consider  
$$f(z_0,z_1, z_2, z_3, z_4)=z_0^3+z_0z_1^2+z_4z_2^{10}+z_2z_3^5+z_3z_4^{14}.$$ The corresponding matrix of exponents $A$ produces the transpose $A^T$ 
\begin{equation*}
    A^T =\begin{bmatrix}
      3 & 1 & 0 & 0 & 0 \\
      0 & 2 & 0 & 0 & 0 \\
      0 & 0 & 10 & 1 & 0 \\
      0 & 0 & 0 & 5 & 1 \\
      0 &  0 & 1 & 0 & 14
    \end{bmatrix}.
\end{equation*}

From Theorems 3.2(2), 3.3(2) and 4.2(1) in \cite{CGL} we conclude that the link $L_{f^T}(\tilde{\mathbf{w}},\tilde{d})$ associated to this matrix satisfies the following properties:
\begin{enumerate}
    \item The link is not well-formed. Moreover, its weight vector is  
    $\tilde{\mathbf{w}}=(701, 2103,342,786,276)$ and the degree is even:  $\tilde{d}=4206.$ Notice that in this case, $2\tilde{w}_2=\tilde{d}.$
    \item This link is a rational homology 7-sphere admitting a Sasaki-Einstein metric.  
\end{enumerate}
The associated polynomial $f^{T}$ to the matrix can be given as the chain-cycle polynomial
$$f^{T}=z_0^3z_1+z_1^2+z_3z_2^{10}+z_4z_3^{5}+z_2z_4^{14}.$$
Now, we consider the link given as a 2-branched cover of $S^9$ with branched locus $L_{f^T}.$ Indeed, it is enough to add one more variable $z_5$ and consider the link 
$$L_g(701, 2103,342,786,276,2103; 4206)=g^{-1}(0)\cap S^{11}$$
associated to the  polynomial 
$$g=z_0^3z_1+z_1^2+z_3z_2^{10}+z_4z_3^{5}+z_2z_4^{14}+z_5^2.$$ From Fact 4(d) we can  alternatively choose the polynomial  $$g_2=z_0^6+z_3z_2^{10}+z_4z_3^{5}+z_2z_4^{14}+z_1^2+z_5^2.$$ 
Since the weight $2\tilde{w}_1=2\tilde{w}_5=\tilde{d},$ we notice that the Sasaki-Reeb cone $\mathfrak{t}^{+}(\mathcal{D}, J)$  of $L_{g_2}$ has dimension 2. 
Here $m_3=701,$ and so it follows from Theorem 4.3, $L_{g_2}$ is diffeomorphic to a Kervaire 9-sphere and of course admits a Sasaki metric of positive Ricci curvature.  
Since (14) is satisfied, we conclude there are no extremal representatives in the whole Sasaki-Reeb cone.
\end{example}

 Links $L_f(\mathbf{w},d)$ such that one of the weights satisfies $2w_i=d$ can be found in the table given in \cite{CGL}. There are 31 elements with that property and with $m_3$ odd. The argument given above applies without changes and hence one obtains 31 new examples of links with Sasaki-Reeb cone of dimension 2 such that the whole cone does not admit extremal Sasaki metrics. These links are  given as 2-fold branched covers of $S^9$  whose branching loci are rational homology 7-spheres arising as  links of  not well-formed orbifolds described by invertible polynomials of type chain-cycle. 
Of course, one can repeat this argument for these 31 elements, that is, one can add more quadratic terms to the original polynomial and apply Theorem 4.3. 
and obtain either rational homology $(4k-1)$-spheres or homotopy $(4k+1)$-spheres both types of manifolds with Sasaki-Reeb cones of dimension $1+\lfloor{\frac{n-3}{2}}\rfloor$  and with no extremal representative in it. All the examples produced here are different Sasakian structures as the ones found in \cite{BvC} since their examples are given as families of Brieskorn-Pham type, actually the numerical data given by the weight vector and the degree of each of the elements of this table forbid singularities of Brieskorn-Pham type. In the following table we present  31  homotopy 9-spheres  (among them, 10 exotic 9-spheres) with no extremal Sasaki metrics with $m_3$ odd.  We label the weights, the polynomial, the degree, $m_3$ and the diffeomorphism type.
\medskip

{\small{    
\begin{center}
\textbf{Table 1: Homotopy 9-spheres that do not admit extremal Sasaki metrics ($m_{3} $ odd)}
\end{center}
 \begin{longtable}{| c | c | c | c | c | } 
 \hline
${\bf {w}}_{Ext}=(w_{0},w_{1},w_{2},w_{3},w_{4},w_{5})$ &  $g_2$  &   $\tilde{d}$   &  $m_{3}$    & Type  \\ \hline  \endfirsthead
\hline
 ${\bf {w}}_{Ext}=(w_{0},w_{1},w_{2},w_{3},w_{4}, w_{5})$ & polynomial $g_2$ & $\tilde{d}$   &  $m_{3}$   & Type \\ \hline  \endhead
 
$(701,2103,342,786,276,2103)$ & $z_{0}^{6}+z_{1}^{2} + z_{3}z_{2}^{10}+z_{4}z_{3}^{5}+z_{2}z_{4}^{14}+z_{5}^{2}$ & $4206$ & $701$  & Kervaire \\ \hline

$(701,2103,246,762,396,2103)$ & $z_{0}^{6}+z_{1}^{2} + z_{3}z_{2}^{14}+z_{4}z_{3}^{5}+z_{2}z_{4}^{10}+z_{5}^{2}$ & $4206$ & $701$ & Kervaire \\ \hline

$(881,2643,1014,216,534,2643)$ & $z_{0}^{6}+z_{1}^{2} + z_{3}z_{2}^{5}+z_{4}z_{3}^{22}+z_{2}z_{4}^{8}+z_{5}^{2}$ & $5286$ & $881$  & Standard \\ \hline

$(881,2643,930,636,198,2643)$ & $z_{0}^{6}+z_{1}^{2} + z_{3}z_{2}^{5}+z_{4}z_{3}^{8}+z_{2}z_{4}^{22}+z_{5}^{2}$ & $5286$ & $881$  & Standard \\ \hline

$(1331,3993,1560,186,918,3993)$ & $z_{0}^{6}+z_{1}^{2} + z_{3}z_{2}^{5}+z_{4}z_{3}^{38}+z_{2}z_{4}^{7}+z_{5}^{2}$ & $7986$ & $1331$ & Kervaire \\ \hline

$(1331,3993,1374,1116,174,3993)$ & $z_{0}^{6}+z_{1}^{2} + z_{3}z_{2}^{5}+z_{4}z_{3}^{7}+z_{2}z_{4}^{38}+z_{5}^{2}$ & $7986$ & $1331$  & Kervaire \\ \hline

$(1457,4371,2112,294,510,4371)$ & $z_{0}^{6}+z_{1}^{2} + z_{3}z_{2}^{4}+z_{4}z_{3}^{28}+z_{2}z_{4}^{13}+z_{5}^{2}$ & $8742$ & $1457$ &   Standard\\ \hline

$(1457,4371,2022,654,240,4371)$ & $z_{0}^{6}+z_{1}^{2} + z_{3}z_{2}^{4}+z_{4}z_{3}^{13}+z_{2}z_{4}^{28}+z_{5}^{2}$ & $8742$ & $1457$ & Standard \\ \hline

$(409,4499,2838,484,770,4499)$ & $z_{0}^{22}+z_{1}^{2} + z_{3}z_{2}^{3}+z_{4}z_{3}^{17}+z_{2}z_{4}^{8}+z_{5}^{2}$ & $8998$ & $409$ &  Standard \\ \hline

$(409,4499,2640,1078,374,4499)$ & $z_{0}^{22}+z_{1}^{2} + z_{3}z_{2}^{3}+z_{4}z_{3}^{8}+z_{2}z_{4}^{17}+z_{5}^{2}$ & $8998$ & $409$ &  Standard \\ \hline

$(2069,6207,3042,246,852,6207)$ & $z_{0}^{6}+z_{1}^{2} + z_{3}z_{2}^{4}+z_{4}z_{3}^{47}+z_{2}z_{4}^{11}+z_{5}^{2}$ & $12414$ & $2069$ &  Kervaire \\ \hline

$(2069,6207,2826,1110,204,6207)$ & $z_{0}^{6}+z_{1}^{2} + z_{3}z_{2}^{4}+z_{4}z_{3}^{11}+z_{2}z_{4}^{47}+z_{5}^{2}$ & $12414$ & $2069$  & Kervaire \\ \hline

$(1297,9079,5950,308,1526,9079)$ & $z_{0}^{14}+z_{1}^{2} + z_{3}z_{2}^{3}+z_{4}z_{3}^{54}+z_{2}z_{4}^{8}+z_{5}^{2}$ & $18158$ & $1297$ &  Standard \\ \hline

$(1297,9079,5306,2240,238,9079)$ & $z_{0}^{14}+z_{1}^{2} + z_{3}z_{2}^{3}+z_{4}z_{3}^{8}+z_{2}z_{4}^{54}+z_{5}^{2}$ & $18158$ & $1297$ &   Standard \\ \hline

$(217,9331,5590,1892,1634,9331)$ & $z_{0}^{86}+z_{1}^{2} + z_{3}z_{2}^{3}+z_{4}z_{3}^{9}+z_{2}z_{4}^{8}+z_{5}^{2}$ & $18662$ & $217$ & Standard \\ \hline

$(217,9331,5504,2150,1462,9331)$ & $z_{0}^{86}+z_{1}^{2} + z_{3}z_{2}^{3}+z_{4}z_{3}^{8}+z_{2}z_{4}^{9}+z_{5}^{2}$ & $18662$ & $217$ &  Standard \\ \hline

$(3401,10203,5046,222,1536,10203)$ & $z_{0}^{6}+z_{1}^{2} + z_{3}z_{2}^{4}+z_{4}z_{3}^{85}+z_{2}z_{4}^{10}+z_{5}^{2}$ & $20406$ & $3401$ &  Standard \\ \hline

$(3401,10203,4596,2022,186,10203)$ & $z_{0}^{6}+z_{1}^{2} + z_{3}z_{2}^{4}+z_{4}z_{3}^{10}+z_{2}z_{4}^{85}+z_{5}^{2}$ & $20406$ & $3401$ &  Standard \\ \hline

$(1135,12485,8206,352,2794,12485)$ & $z_{0}^{22}+z_{1}^{2} + z_{3}z_{2}^{3}+z_{4}z_{3}^{63}+z_{2}z_{4}^{6}+z_{5}^{2}$ & $24970$ & $1135$ &  Standard \\ \hline

$(1135,12485,6952,4114,286,12485)$ & $z_{0}^{22}+z_{1}^{2} + z_{3}z_{2}^{3}+z_{4}z_{3}^{6}+z_{2}z_{4}^{63}+z_{5}^{2}$ & $24970$ & $1135$ &  Standard \\ \hline

$(1505,13545,5094,6714,234,13545)$ & $z_{0}^{18}+z_{1}^{2} + z_{3}z_{2}^{4}+z_{4}z_{3}^{4}+z_{2}z_{4}^{94}+z_{5}^{2}$ & $27090$ & $1505$ &  Standard\\ \hline

$(141,19599,8618,4726,6116,19599)$ & $z_{0}^{278}+z_{1}^{2} + z_{3}z_{2}^{4}+z_{4}z_{3}^{7}+z_{2}z_{4}^{5}+z_{5}^{2}$ & $39198$ & $141$  & Kervaire \\ \hline

$(141,19599,8618,4726,6116,19599)$ & $z_{0}^{278}+z_{1}^{2} + z_{3}z_{2}^{4}+z_{4}z_{3}^{5}+z_{2}z_{4}^{7}+z_{5}^{2}$ & $39198$ & $141$  & Kervaire \\ \hline

$(1351,25669,16948,494,6878,25669)$ & $z_{0}^{38}+z_{1}^{2} + z_{3}z_{2}^{3}+z_{4}z_{3}^{90}+z_{2}z_{4}^{5}+z_{5}^{2}$ & $51338$ & $1351$  & Standard \\ \hline

$(1351,25669,13718,10184,418,25669)$ & $z_{0}^{38}+z_{1}^{2} + z_{3}z_{2}^{3}+z_{4}z_{3}^{5}+z_{2}z_{4}^{90}+z_{5}^{2}$ & $51338$ & $1351$  & Standard \\ \hline

$(177,30975,11900,14350,4550,30975)$ & $z_{0}^{350}+z_{1}^{2} + z_{3}z_{2}^{4}+z_{4}z_{3}^{4}+z_{2}z_{4}^{11}+z_{5}^{2}$ & $61950$ & $177$  & Standard \\ \hline

$(193,36863,21774,8404,6494,36863)$ & $z_{0}^{382}+z_{1}^{2} + z_{3}z_{2}^{3}+z_{4}z_{3}^{8}+z_{2}z_{4}^{8}+z_{5}^{2}$ & $73726$ & $193$ &  Standard \\ \hline

$(217,46655,28810,6880,10750,46655)$ & $z_{0}^{430}+z_{1}^{2} + z_{3}z_{2}^{3}+z_{4}z_{3}^{12}+z_{2}z_{4}^{6}+z_{5}^{2}$ & $93310$ & $217$ &  Standard \\ \hline

$(217,46655,26230,14620,5590,46655)$ & $z_{0}^{430}+z_{1}^{2} + z_{3}z_{2}^{3}+z_{4}z_{3}^{6}+z_{2}z_{4}^{12}+z_{5}^{2}$ & $93310$ & $217$ &  Standard \\ \hline

$(301,89999,57408,7774,24518,89999)$ & $z_{0}^{598}+z_{1}^{2} + z_{3}z_{2}^{3}+z_{4}z_{3}^{20}+z_{2}z_{4}^{5}+z_{5}^{2}$ & $179998$ & $301$ &  Kervaire \\ \hline

$(301,89999,48438,34684,6578,89999)$ & $z_{0}^{598}+z_{1}^{2} + z_{3}z_{2}^{3}+z_{4}z_{3}^{5}+z_{2}z_{4}^{20}+z_{5}^{2}$ & $179998$ & $301$ &  Kervaire \\\hline
\end{longtable}
}}
\medskip

In case $m_3$ is even, we have the following result: 

\begin{thm} Consider a polynomial of chain-cycle type $$f=z_{0}^{a_{0}}+z_{0}z_{1}^{a_{1}}+z_{4}z_{2}^{a_{2}}+z_{2}z_{3}^{a_{3}}+z_{3}z_{4}^{a_{4}}$$ with $a_1=2$ whose link $L_f$ is a rational homology sphere and  that cuts out a projective hypersurface in $\mathbb{P}(w_0, w_1, \ldots , w_4)$ such that $$(w_{0},w_{1},w_{2},w_{3},w_{4})=(m_{3}v_{0},m_{3}v_{1},m_{2}v_{2},m_{2}v_{3},m_{2}v_{4}),$$ with  $\gcd(m_{2},m_{3})=1$ and $d=m_3m_2$ and $m_{3}$ even. Then 
the polynomial $$g=f^T+z_5^2+z_6^2$$ with $f^T$ the Berglund-Hübsch dual of $f,$  determines a link $L_g$  which is of the form $S^{4}\times S^{5}$. In particular, all the elements on the list in the table in the appendix in \cite{CGL} satisfying the conditions given above do not admit extremal Sasaki metrics in the whole Sasaki-Reeb cone which is of dimension 2. 
    \end{thm}

    \begin{pf}
        In this case, we can simplify the weights $\tilde{w}_{i}'s$:
        $$\tilde{w}_{0}=\dfrac{m_{3}}{2}, \ \ \tilde{w}_{1}=\dfrac{m_{3}m_{2}}{2}, \ \ \tilde{w}_{2}=m_{2}\tilde{v}_{2}, \ \ \tilde{w}_{3}=m_{2}\tilde{v}_{3}, \ \ \tilde{w}_{4}=m_{2}\tilde{v}_{4}, \ \ \mbox{ and } \ \ \tilde{w}_{5}=\dfrac{m_{3}m_{2}}{2}.$$
        Also, we have the degree $\tilde{d}=m_{3}m_{2}$ and the index $\tilde{I}=1$.
        Next, we calculate the numbers $u_{i}=\dfrac{\tilde{d}}{\gcd(\tilde{d},\tilde{w}_{i})}$ and $v_{i}=\dfrac{\tilde{w}_{i}}{\gcd(\tilde{d},\tilde{w}_{i})}$ for these weights:
\begin{center}
      \begin{tabular}{| c || c | c | c| c | c | c |}
\hline
$i$ & $0$ & $1$ & $2$ & $3$ & $4$ & $5$\\ \hline
$u_{i}$ &  $2m_{2}$ & $2$ & $m_{3}$ & $m_{3}$ & $m_{3}$ & $2$ \\ \hline
$v_{i}$ & $1$ & $1$ & $\tilde{v}_{2}$ & $\tilde{v}_{3}$ & $\tilde{v}_{4}$ & $1$ 
\\ \hline
\end{tabular}
\end{center}
The Alexander polynomial for $g$ is:

$$\mbox{div } \Delta_{g} = \prod_{i=0}^{5}\left( \dfrac{\Lambda_{u_{i}}}{v_{i}}-\Lambda_{1} \right).$$
Using the relation $\Lambda_{a}\Lambda_{b}=\gcd(a,b)\Lambda_{\mbox{lcm}(a,b)}$, we reduce:
$$\mbox{div } \Delta_{g} = (\Lambda_{2m_{2}}-\Lambda_{1})(\Lambda_{m_{3}}-\Lambda_{1}).$$
As $m_{3}$ is even, we obtain
$$\mbox{div } \Delta_{g} = 2\Lambda_{m_{2}m_{3}}-\Lambda_{m_{3}}-\Lambda_{2m_{2}}+\Lambda_{1}.$$
Therefore, the Alexander polynomial for $g$ is 
$$\Delta_{g}(t)=\dfrac{(t^{m_{2}m_{3}}-1)^{2}(t-1)}{(t^{m_{3}}-1)(t^{2m_{2}}-1)}.$$
Since the expression above only has a factor $(t-1)$, then $b_{4}(L_{g})=1$.

Next, we calculate the torsion for $L_{g}$. By the Orlik's algorithm, we obtain, $c_{\emptyset}=2, c_{015}=\dfrac{m_{3}}{2}$, $c_{12345}=m_{2}$, and $c_{i_{1}\dots i_{s}}=1$ in other cases. Since $k_{\emptyset}=0$, we just need to calculate $k_{015}$ and $k_{12345}$. Using the formula for $k_{i_{1}\dots i_{s}}$, we obtain $k_{015}=0$ and $k_{12345}=0$. Therefore, $L_{g}$ has the form $S^{4}\times S^{5}$.
\hfill$\square$
\end{pf}

In the following table, we present all links that do not admit extremal Sasaki metrics from the list given in the table in the appendix in \cite{CGL} where $m_{3}$ is even and $a_{1}=2$. Here the dimension of the Sasaki-Reeb cone equals two. The inequality (13) is easily verified for these polynomials.
\medskip

\begin{center}
\noindent{\bf Table 2: Links of the form $S^{4}\times S^{5}$  with no extremal Sasaki metrics ($m_3$ even)} 
\end{center}
\begin{longtable}{| c | c | c | c | }
 \hline
$\mathbf{w}_{Ext}=(w_{0},w_{1},w_{2},w_{3},w_{4},w_{5})$ & $g_2$ &   $\tilde{d}$   &  $\Delta_{g}(-1)=m_{3}$   \\ 
\hline \hline 
\endfirsthead
\hline
 $\mathbf{w}=(w_{0},w_{1},w_{2},w_{3},w_{4})$ & polynomial $g=f^{T}+z_{5}^2$ & $\tilde{d}$   &  $\Delta_{g}(-1)=m_{3}$  \\ \hline \hline \endhead
$(1766,8830,5835,155,1075,8830)$ & $z_{0}^{10}+z_{1}^{2} + z_{3}z_{2}^{3}+z_{4}z_{3}^{107}+z_{2}z_{4}^{11}+z_{5}^{2}$ & $17660$ & $3532$  \\ \hline

$(1766,8830,5355,1595,115,8830)$ & $z_{0}^{10}+z_{1}^{2} + z_{3}z_{2}^{3}+z_{4}z_{3}^{11}+z_{2}z_{4}^{107}+z_{5}^{2}$ & $17660$ & $3532$   \\ \hline

$(1208,20536,13617,221,5491,20536)$ & $z_{0}^{34}+z_{1}^{2} + z_{3}z_{2}^{3}+z_{4}z_{3}^{161}+z_{2}z_{4}^{5}+z_{5}^{2}$ & $41072$ & $2416$  \\ \hline

$(1208,20536,10965,8177,187,20536)$ & $z_{0}^{34}+z_{1}^{2} + z_{3}z_{2}^{3}+z_{4}z_{3}^{5}+z_{2}z_{4}^{161}+z_{5}^{2}$ & $41072$ & $2416$ \\ \hline

$(158,24806,15857,2041,6751,24806)$ & $z_{0}^{314}+z_{1}^{2} + z_{3}z_{2}^{3}+z_{4}z_{3}^{21}+z_{2}z_{4}^{5}+z_{5}^{2}$ & $49612$ & $316$  \\ \hline

$(158,24806,13345,9577,1727,24806)$ & $z_{0}^{314}+z_{1}^{2} + z_{3}z_{2}^{3}+z_{4}z_{3}^{5}+z_{2}z_{4}^{21}+z_{5}^{2}$ & $49612$ & $316$ \\ \hline

\end{longtable}

\medskip

\subsection{Examples derived from Thom-Sebastiani sums of chain type singularities}
In this subsection, we found examples of Sasakian manifolds where the Lichnerowicz inequality is preserved via the Berglund-H\"ubsch rule, that is, we produce links such that if the link $L_f$ satisfies the Lichnerowicz obstruction  then its Berglund-H\"ubsch dual $L_{f^T}$ preserves it as well, and additionally, $L_{f^T}$ always admits a  deformation $L_{f_{BP}^T}$ in its local moduli which is obstructed to admitting  extremal Sasaki metrics in its whole cone. All these examples  are either  homotopy spheres or rational homology spheres. 

Let us begin considering polynomials of the form BP-Chain  
$$
f=z_{0}^{a_{0}}+z_{1}^{a_{1}}+z_{1}z_{2}^{2}+z_{3}^{a_{3}}+z_{3}z_{4}^{2}
$$
and add some arithmetic conditions on the exponents to obtain links that are rational homology 7-spheres

\begin{lemma}
Given the invertible polynomial of the form 
\begin{equation}
f=z_{0}^{a_{0}}+z_{1}^{a_{1}}+z_{1}z_{2}^{2}+z_{3}^{a_{3}}+z_{3}z_{4}^{2}
\end{equation}
with associated weight vector ${\bf{w}}=(w_{0},w_{1},w_{2},w_{3},w_{4})$.
If $a_{0}, a_{1}$ and $a_{3}$ are pairwise  relatively prime integers and $a_{0}$ is odd, then the link $L_{f^{T}}$ associated to the weight vector $\tilde{\bf{w}}$ obtained by the Berglund-Hübsch transpose rule is a rational homology 7-sphere with $ H_{3}(L_{f^{T}},\mathbb{Z}) =\mathbb{Z}_{a_{0}}$ and Milnor number $\mu=(a_{0}-1)(2a_{1}-1)(2a_{3}-1).$ 
\end{lemma}
\begin{pf}
    Let  $\textbf{w}=(w_{0},w_{1},w_{2},w_{3},w_{4})$ and $d$ the weight vector and degree associated to $f$. Its exponent matrix is given by 
    \begin{equation*}
       A= \begin{bmatrix}
          a_{0} & 0 & 0 & 0 & 0  \\
      0 & a_{1} & 0 & 0 & 0 \\
      0 & 1 & 2 & 0 & 0 \\
      0 & 0 & 0 & a_{3} & 0 \\
      0 &  0 & 0 & 1 & 2   
        \end{bmatrix}.
        \end{equation*}
 Applying the Berglund-Hübsch transpose rule, we obtain
$$    A^T= \begin{bmatrix}
          a_{0} & 0 & 0 & 0 & 0  \\
      0 & a_{1} & 1 & 0 & 0 \\
      0 & 0 & 2 & 0 & 0 \\
      0 & 0 & 0 & a_{3} & 1 \\
      0 &  0 & 0 & 0 & 2   
        \end{bmatrix}.
        $$
Now we are looking for weights $\tilde{\textbf{w}}=(\tilde{w}_{0},\tilde{w}_{1},\tilde{w}_{2},\tilde{w}_{3},\tilde{w}_{4})$ and degree $\tilde{d}$ that satisfy the matrix equation $A^T\tilde{\bf w}^{T}=\tilde{D}^{T}$, where $\tilde{D}=(\tilde{d},\tilde{d},\tilde{d},\tilde{d},\tilde{d})$. Notice that this is equivalent to solving the equation:
\begin{equation}
        \begin{bmatrix}
          a_{0} & 0 & 0 & 0 & 0 & -1 \\
      0 & a_{1} & 1 & 0 & 0 & -1\\
      0 & 0 & 2 & 0 & 0 & -1\\
      0 & 0 & 0 & a_{3} & 1 & -1\\
      0 &  0 & 0 & 0 & 2 & -1  
        \end{bmatrix} \begin{bmatrix}
            \Tilde{w}_{0} \\
            \Tilde{w}_{1} \\
            \Tilde{w}_{2} \\
            \Tilde{w}_{3} \\
            \Tilde{w}_{4} \\
            \Tilde{d}
        \end{bmatrix} = \begin{bmatrix}
            0 \\ 
            0 \\
            0 \\
            0 \\
            0
        \end{bmatrix}.
    \end{equation}
    We suggest the following non-trivial solution with integer entries
    \begin{equation}
        \begin{bmatrix}
            \tilde{w}_{0} & \tilde{w}_{1} & \tilde{w}_{2} & \tilde{w}_{3} & \tilde{w}_{4} & \tilde{d} 
        \end{bmatrix}^{T} = \begin{bmatrix}
            2w_{0} & w_{1} & d & w_{3} & d & 2d 
        \end{bmatrix}^{T}
    \end{equation}

    Using the weights proposed above, we can determine the numbers $u_{i}$ and $v_{i}$, for $i=0,1,\dots,4$. First, notice  that the values of $u_{i}$ and $v_{i}$ are well defined. Indeed, since $A^T$ is an invertible matrix, the space of solutions of Equation (20) are 1-dimensional. Therefore, any solution with integer entries of this system of equations is a multiple of the column matrix given in equality (21):
    \begin{equation*}
        \begin{bmatrix}
            \tilde{w}_{0}' & \tilde{w}_{1}' & \tilde{w}_{2}' & \tilde{w}_{3}' & \tilde{w}_{4}' & \tilde{d}' 
        \end{bmatrix}^{T} = \lambda\begin{bmatrix}
            \tilde{w}_{0} & \tilde{w}_{1} & \tilde{w}_{2} & \tilde{w}_{3} & \tilde{w}_{4} & \tilde{d} 
        \end{bmatrix}^{T}
    \end{equation*}
    Next we compute the values $u_{i}'$ and $v_{i}'$ with respect to the weights $\tilde{w}_{i}'$ and degree $\tilde{d}'$. Thus, we obtain
    {\small{
    $$u_{i}'=\dfrac{\tilde{d}'}{\gcd(\tilde{d}',\tilde{w}_{i}')}=\dfrac{\lambda\tilde{d}}{\lambda\gcd(\tilde{d},\tilde{w}_{i})}=u_{i} \ \ \mbox{ and } \ \ v_{i}'=\dfrac{\tilde{w}_{i}'}{\gcd(\tilde{d}',\tilde{w}_{i}')}=\dfrac{\lambda\tilde{w}_{i}}{\lambda\gcd(\tilde{d},\tilde{w}_{i})}=v_{i}.$$}}
    Hence, it suffices to work with the weights $\tilde{w}_{i}'s$ and the degree $\tilde{d}$ obtained in equality (21). 
    Then, we calculate the values $u_{i}$ and $v_{i}$:
    \begin{center}
      \begin{tabular}{| c || c | c | c| c | c | }
\hline
$i$ & $0$ & $1$ & $2$ & $3$ & $4$ \\ \hline
$u_{i}$ &  $a_{0}$ & $2a_{1}$ & $2$ & $2a_{3}$ & $2$  \\ \hline
$v_{i}$ & $1$ & $1$ & $1$ & $1$ & $1$  
\\ \hline
\end{tabular}.
\end{center}

From $$\mbox{div } \Delta_{f^{T}} = \prod_{i=0}^{4}\left( \dfrac{\Lambda_{u_{i}}}{v_{i}}-\Lambda_{1} \right)$$ we obtain the following equality

\begin{equation}
\mbox{div } \Delta_{f^{T}} = (\Lambda_{a_{0}}-\Lambda_{1})(\Lambda_{2a_{1}}-\Lambda_{1})(\Lambda_{2}-\Lambda_{1})(\Lambda_{2a_{3}}-\Lambda_{1})(\Lambda_{2}-\Lambda_{1})
\end{equation}

Since $a_{0},a_{1}$ and $a_{3}$ are co-prime, $a_{0}$ is odd and $(\Lambda_{2}-\Lambda_{1})(\Lambda_{2}-\Lambda_{1})=1$, Equation (22) reduces  to
\begin{align*}
\mbox{div } \Delta_{f^{T}} & = (\Lambda_{a_{0}}-\Lambda_{1})(\Lambda_{2a_{1}}-\Lambda_{1})(\Lambda_{2a_{3}}-\Lambda_{1}) \\
& = 2\Lambda_{2a_{0}a_{1}a_{3}}+\Lambda_{a_{0}}+\Lambda_{2a_{1}}+\Lambda_{2a_{3}}-\Lambda_{2a_{0}a_{1}}-\Lambda_{2a_{0}a_{3}}-\Lambda_{2a_{1}a_{3}}-\Lambda_{1}.
\end{align*}
Therefore, the Alexander polynomial for $f^{T}$ is 
$$\Delta_{f^{T}}(t)=\dfrac{(t^{2a_{0}a_{1}a_{3}}-1)^{2}(t^{a_{0}}-1)(t^{2a_{1}}-1)(t^{2a_{3}}-1)}{(t^{2a_{0}a_{1}}-1)(t^{2a_{0}a_{3}}-1)(t^{2a_{1}a_{3}}-1)^{2}(t-1)}.$$
Since $\Delta_{f^{T}}(1)=a_{0}\neq 1$, the link $L_{f^{T}}$ is a rational homology sphere. Now, we calculate the torsion of this link. Using Orlik's formula, we get the quantities $c_{\emptyset}=1, c_{0}=2, c_{0124}=a_{3}, c_{0234}=a_{1}, c_{1234}=a_{0}$ and $c_{i_{1}\dots i_{s}}=1$ in other cases. Since $k_{0}=0$, we only need to calculate the following terms: $k_{0124}=0, k_{0234}=0$ and $k_{1234}=1$. Thus  $H_{3}(L_{f^{T}},\mathbb{Z})=\mathbb{Z}_{a_{0}}$.

Finally, we have that the Milnor number is given by
\begin{align*}
  \mu(L_{f^{T}}) & = \left( \dfrac{\tilde{d}-\tilde{w}_{0}}{\tilde{w}_{0}}\right)\left( \dfrac{\tilde{d}-\tilde{w}_{1}}{\tilde{w}_{1}}\right) \left( \dfrac{\tilde{d}-\tilde{w}_{2}}{\tilde{w}_{2}}\right) \left( \dfrac{\tilde{d}-\tilde{w}_{3}}{\tilde{w}_{3}}\right) \left( \dfrac{\tilde{d}-\tilde{w}_{4}}{\tilde{w}_{4}}\right) \\
  & = (a_{0}-1)\left(\dfrac{2d-w_{1}}{w_{1}}\right)\left(\dfrac{2d-d}{d}\right)\left(\dfrac{2d-w_{3}}{w_{3}}\right)\left(\dfrac{2d-d}{d}\right) \\
  & = (a_{0}-1)(2a_{1}-1)(2a_{3}-1).
\end{align*}
\hfill$\square$
\end{pf}

\begin{rmk}
    From a similar argument to the one given in Remark 3.1 one can obtain for the weight vector $\bf \tilde{w}$ obtained in the proof of the previous lemma  a polynomial 
of Brieskorn-Pham type $$f_{BP}^{T}=z_{0}^{a_{0}}+z_{1}^{2a_{1}}+z_{2}^{2}+z_{3}^{2a_{3}}+z_{4}^{2}\in H^0(\mathbb{P}({\bf\tilde{w}}), \mathcal{O}(\tilde{d})),$$ with degree $\tilde{d}=\deg(f^T)=2d$ with $f$ given as (19) with no assumption on the exponents of $f.$
\end{rmk}

Next, we will show that under certain conditions on the weights of the link $L_{f}$,  the rational homology spheres $L_{f_{BP}^{T}}$ do not admit extremal Sasaki metrics. 


\begin{thm}
    Consider the invertible polynomial $$f=z_{0}^{a_{0}}+z_{1}^{a_{1}}+z_{1}z_{2}^{2}+z_{3}^{a_{3}}+z_{3}z_{4}^{2}$$ 
    with $a_1\geq 3$ or $a_3\geq 3$ and $a_0$ odd. If the corresponding link $L_f$ is such that the Lichnerowicz obstruction is satisfied, then the Berglund-H\"ubsch dual 
    $L_{f^T}$ of $L_f$ satisfies this obstruction as well. Moreover, $L_{f^{T}_{BP}}$,  a perturbation of  $L_{f^T}$, which is defined via $$f_{BP}^{T}=z_{0}^{a_{0}}+z_{1}^{2a_{1}}+z_{2}^{2}+z_{3}^{2a_{3}}+z_{4}^{2}$$ given in Remark 4.2 satisfies the generalized Lichnerowicz obstruction (13) and therefore the link $L_{f_{BP}^{T}}$ does not admit extremal Sasaki metrics in its whole Sasaki-Reeb cone which is of dimension 2. Furthermore, if additionally $a_0, a_1$ and $a_3$ are coprime, then $L_{f_{BP}^{T}}$ is a rational homology 7-sphere with no extremal metrics in its whole Sasaki-Reeb cone which is of dimension 2.
    \end{thm}
    \begin{pf}
     For the given $f$, we can write its weights $w_{i}$'s as
       {\small{
       \begin{equation}
          w_{0}=\dfrac{d}{a_{0}}, \ \ w_{1}=\dfrac{d}{a_{1}}, \ \ w_{2}=\dfrac{(a_{1}-1)d}{2a_{1}}, \ \ w_{3}=\dfrac{d}{a_{3}} \ \ \mbox{ and } \ \ w_{4}=\dfrac{(a_{3}-1)d}{2a_{3}}  
       \end{equation}
       }}
       where the index $I$ can be obtained as
       $$I=|{\bf w}|-d=w_{0}+\dfrac{w_{1}}{2}+\dfrac{w_{3}}{2}.$$
       Since $a_{1}\geq3$ or $a_{3}\geq3$, we have  that $$w_{2}=\dfrac{(a_{1}-1)d}{2a_{1}}\geq \dfrac{d}{a_{1}}=w_{1} \ \ \mbox{ or } \ \ w_{4}=\dfrac{(a_{3}-1)d}{2a_{3}}\geq \dfrac{d}{a_{3}}=w_{3}.$$
       Therefore, we have 
       \begin{equation}
           \min_{i}w_{i} \in \{ w_{0},w_{1},w_{3}\}.
       \end{equation}
        On the other hand, from Equation (21), we obtain the weights of $\bf \tilde{w}$ and the degree $\tilde{d}$ for the polynomial $f^{T}$:
        {\small{
        \begin{equation}
            \tilde{w}_{0}=2w_{0}=\dfrac{2d}{a_{0}}, \ \ \tilde{w}_{1}=w_{1}=\dfrac{d}{a_{1}}, \ \ \tilde{w}_{2}=d, \ \ \tilde{w}_{3}=w_{3}=\dfrac{d}{a_{3}}, \ \ \tilde{w}_{4}=d \ \ \mbox{ and } \ \ \tilde{d}=2d
        \end{equation}
        }}
       where the index $\tilde{I}$ of $L_{f^{T}}$ can be calculated as
       $$\tilde{I}=| \tilde{\mathbf{w}} | - \tilde{d} = \tilde{w}_{0}+\tilde{w}_{1}+\tilde{w}_{3}=2w_{0}+w_{1}+w_{3}.$$
       From equality (25) we have 
       \begin{equation}
           \min_{i}\tilde{w}_{i}\in \{\tilde{w}_{0},\tilde{w}_{1},\tilde{w}_{3}\}.
       \end{equation}
       Here we consider the following cases:
       \begin{itemize}
           \item[a)] We suppose that $\min_{i}w_{i} = w_{1}$ and the link $L_{f}$ verifies the Lichnerowicz inequality:
           $$I=w_{0}+\dfrac{w_{1}}{2}+\dfrac{w_{3}}{2}>4w_{1}.$$
           Using equality (25)  we have
           $$\tilde{w}_{0}+\tilde{w}_{1}+\tilde{w}_{3}> w_{0}+\dfrac{w_{1}}{2}+\dfrac{w_{3}}{2}>4w_{1} = 4\tilde{w}_{1}.$$ 
           Then $\tilde{I}>4 \tilde{w}_{1}\geq 4 \min \{\tilde{w}_{i}\}.$ Thus $L_{f^T}$ satisfies the Lichnerowicz obstruction. 
           \item[b)] We suppose  that $\min_{i}w_{i} = w_{3}$ and the link $L_{f}$ verifies the Lichnerowicz inequality. The proof of this case is similar to case a).
           \item[c)] We suppose that $\min_{i}w_{i} = w_{0}$ and the link $L_{f}$ satisfies  the Lichnerowicz inequality. Then we have
           $$I=w_{0}+\dfrac{w_{1}}{2}+\dfrac{w_{3}}{2}>4w_{0}.$$ 
           Using equality (25) and multiplying by $2$ the previous inequality, we obtain 
           $$\tilde{w}_{0}+\tilde{w}_{1}+\tilde{w}_{3}=2w_{0}+w_{1}+w_{3}>4(2w_{0})=4\tilde{w}_{0}\geq 4\min\{\tilde{w}_{i}\}.$$ 
           \end{itemize}

By Remark 4.2, the weight vector  $\tilde{\bf w}$ admits a Brieskorn-Pham polynomial of the form        $$f_{BP}^{T}=z_{0}^{a_{0}}+z_{1}^{2a_{1}}+z_{2}^{2}+z_{3}^{2a_{3}}+z_{4}^{2},$$   and notice that Lichnerowicz obstruction  satisfied by $L_{f^T}$ forces  the generalized Lichnerowicz inequality (13) in $f_{BP}^{T}.$ Thus, $L_{f_{BP}^{T}}$ does not admit extremal Sasaki metrics. The last claim follows from Lemma 4.5
\hfill$\square$
 \end{pf}

\begin{rmk}
    When we add an even number of quadratic terms $z_{i}^{2}$ to the polynomial $f_{BP}^{T}$, we obtain the polynomial
    $$f_{Ext}^{T}=f_{BP}^{T}+z_{5}^{2}+\dots + z_{2n}^{2}.$$ The dimension of the Sasaki-Reeb cone of $L_{f_{Ext}^{T}}$ increases by one for each two terms added, so we obtain a Sasaki-Reeb cone of dimension $n$.  Since $L_{f_{BP}^{T}}$ is a rational homology sphere, it is clear that  $L_{f_{Ext}^{T}}$ is also a rational homology sphere (see Corollary 9.5.3 in \cite{BBG}). Moreover,
    the link $L_{f_{Ext}^{T}}$ does not admit extremal Sasaki metrics in its Sasaki-Reeb cone since $L_{f_{BP}^{T}}$ does not admit extremal Sasaki metrics either. Indeed, we can apply the formula (13) to $f_{Ext}^{T}$. In our case, we always have $k=2$, thus we must prove
    \begin{equation}
        \tilde{w}_{0}+\tilde{w}_{1}+\tilde{w}_{3}+(2n-4)\dfrac{\tilde{d}}{2}\geq (2n)\min_{i}\tilde{w}_{i}.
    \end{equation}
    Since $\tilde{w}_{i}=\dfrac{d}{2}$ for all $i\geq5$, we have that $\min_{i}\tilde{w}_{i}\in\{ \tilde{w}_{0},\tilde{w}_{1},\tilde{w}_{3}\}$. Thus, we obtain $$\min_{i}\tilde{w}_{i}<\dfrac{\tilde{d}}{2}.$$ By the previous lemma, we also have $\tilde{w}_{0}+\tilde{w}_{1}+\tilde{w}_{3}>4\min_{i}\tilde{w}_{i}$. Putting together these results, we obtain
    \begin{align*}
        \tilde{w}_{0}+\tilde{w}_{1}+\tilde{w}_{3}+(2n-4)\dfrac{\tilde{d}}{2} & > 4\min_{i}\tilde{w}_{i}+(2n-4)\min_{i}\tilde{w}_{i} \\
        & = (2n)\min_{i}\tilde{w}_{i}.
    \end{align*}
    Therefore, the inequality (27) holds. 
\end{rmk}

Now, let us consider examples with the topology of a homotopy 7-sphere.

\begin{lemma}
Let $g$ be the invertible polynomial of the form 
\begin{equation}
  g= z_{0}^{2}+z_{1}^{a_{1}}+z_{2}^{a_{2}}+z_{3}^{a_{3}}+z_{3}z_{4}^{2}  
\end{equation}
with associated weight vector ${\bf{w}}=(w_{0},w_{1},w_{2},w_{3},w_{4})$.
If $a_{1},a_{2}$ and $a_{3}$ are co-prime and $a_{1},a_{2}$ are odd numbers, then the link $L_{g^{T}}$ associated to the weight vector $\tilde{\bf{w}}$ obtained by the Berglund-Hübsch transpose rule is a homotopy sphere. Moreover, its Milnor number $\mu(L_{g^{T}})$ equals $(a_{1}-1)(a_{2}-1)(2a_{3}-1)$.
\end{lemma}
\begin{pf}
    The  matrix of exponents associated to polynomial $g$ defined in (28) is given by
    \begin{equation}
       A_{g}= \begin{bmatrix}
          2 & 0 & 0 & 0 & 0  \\
      0 & a_{1} & 0 & 0 & 0 \\
      0 & 0 & a_{2} & 0 & 0 \\
      0 & 0 & 0 & a_{3} & 0 \\
      0 &  0 & 0 & 1 & 2   
        \end{bmatrix}.
        \end{equation}
Applying the Berglund-Hübsch transpose rule, we have the exponent matrix
\begin{equation*}
       A_{g^{T}}= \begin{bmatrix}
          2 & 0 & 0 & 0 & 0  \\
      0 & a_{1} & 0 & 0 & 0 \\
      0 & 0 & a_{2} & 0 & 0 \\
      0 & 0 & 0 & a_{3} & 1 \\
      0 &  0 & 0 & 0 & 2   
        \end{bmatrix}.
        \end{equation*}
The dual $g^{T}$ is given by
\begin{equation*}
  g^{T}= z_{0}^{2}+z_{1}^{a_{1}}+z_{2}^{a_{2}}+z_{3}^{a_{3}}z_{4}+z_{4}^{2}.  
\end{equation*}
Now, we will find weights $\tilde{\bf w}=(\tilde{w}_{0},\tilde{w}_{1},\tilde{w}_{2},\tilde{w}_{3},\tilde{w}_{4})$ and degree $\tilde{d}$ that solve the matrix equation $A_{g^{T}}\tilde{w}^{T}=\tilde{D}^{T}$, where $\tilde{D}=(\tilde{d},\tilde{d},\tilde{d},\tilde{d},\tilde{d})$. In a similar way as we have done above, these weights are obtained solving the following equation
\begin{equation}
        \begin{bmatrix}
          2 & 0 & 0 & 0 & 0 & -1 \\
      0 & a_{1} & 0 & 0 & 0 & -1\\
      0 & 0 & a_{2} & 0 & 0 & -1\\
      0 & 0 & 0 & a_{3} & 1 & -1\\
      0 &  0 & 0 & 0 & 2 & -1  
        \end{bmatrix} \begin{bmatrix}
            \Tilde{w}_{0} \\
            \Tilde{w}_{1} \\
            \Tilde{w}_{2} \\
            \Tilde{w}_{3} \\
            \Tilde{w}_{4} \\
            \Tilde{d}
        \end{bmatrix} = \begin{bmatrix}
            0 \\ 
            0 \\
            0 \\
            0 \\
            0
        \end{bmatrix}.
    \end{equation}
Then we propose a solution with integer entries for Equation (30), which is given by
\begin{equation}
    \tilde{w}_{0}=d, \ \ \tilde{w}_{1}=2w_{1}, \ \ \tilde{w}_{2}=2w_{2}, \ \ \tilde{w}_{3}=w_{3}, \ \ \tilde{w}_{4}=d \ \ \mbox{ and } \ \ \tilde{d}=2d.
\end{equation}
As we have mentioned above, we can find  integer weights $\tilde{w}_{i}$'s and degree $\tilde{d}$ that solve Equation (30); nevertheless, to calculate the homology of $L_{g^{T}}$ it is sufficient to work with the weights and degree given in equality (31). Thus, we obtain the numbers $u_{i}$ and $v_{i}$:
\begin{center}
      \begin{tabular}{| c || c | c | c| c | c | }
\hline
$i$ & $0$ & $1$ & $2$ & $3$ & $4$ \\ \hline
$u_{i}$ &  $2$ & $a_{1}$ & $a_{2}$ & $2a_{3}$ & $2$  \\ \hline
$v_{i}$ & $1$ & $1$ & $1$ & $1$ & $1$  
\\ \hline
\end{tabular}.
\end{center}
Then we calculate the Alexander polynomial of $g^{T}$:
$$\mbox{div } \Delta_{g^{T}} = \prod_{i=0}^{4}\left( \dfrac{\Lambda_{u_{i}}}{v_{i}}-\Lambda_{1} \right)= (\Lambda_{2}-\Lambda_{1})(\Lambda_{a_{1}}-\Lambda_{1})(\Lambda_{a_{2}}-\Lambda_{1})(\Lambda_{2a_{3}}-\Lambda_{1})(\Lambda_{2}-\Lambda_{1}).$$

Since $a_{1},a_{2}$ and $a_{3}$ are co-primes, where $a_{1}$ and $a_{2}$ are odd, and $(\Lambda_{2}-\Lambda_{1})(\Lambda_{2}-\Lambda_{1})=1$, the expression above can be written as 
\begin{align*}
\mbox{div } \Delta_{g^{T}} & = (\Lambda_{a_{1}}-\Lambda_{1})(\Lambda_{a_{2}}-\Lambda_{1})(\Lambda_{2a_{3}}-\Lambda_{1}) \\
& = \Lambda_{2a_{1}a_{2}a_{3}}+\Lambda_{a_{1}}+\Lambda_{a_{2}}+\Lambda_{2a_{3}}-\Lambda_{a_{1}a_{2}}-\Lambda_{2a_{1}a_{3}}-\Lambda_{2a_{2}a_{3}}-\Lambda_{1}.
\end{align*}
Then the Alexander polynomial of $g^{T}$ is given by
$$\Delta_{g^{T}}(t)=\dfrac{(t^{2a_{1}a_{2}a_{3}}-1)(t^{a_{1}}-1)(t^{a_{2}}-1)(t^{2a_{3}}-1)}{(t^{a_{1}a_{2}}-1)(t^{2a_{1}a_{3}}-1)(t^{2a_{2}a_{3}}-1)(t-1)}.$$
Evaluating at $t=1$:
$$ \Delta_{g^{T}}(1)=\dfrac{(2a_{1}a_{2}a_{3})(a_{1})(a_{2})(2a_{3})}{(a_{1}a_{2})(2a_{1}a_{3})(2a_{2}a_{3})}=1.$$
Then, the link $L_{g^{T}}$ is an homotopy 7-sphere. Finally, we calculate the Milnor number for this link
\begin{align*}
  \mu(L_{g^{T}}) & = \left( \dfrac{\tilde{d}-\tilde{w}_{0}}{\tilde{w}_{0}}\right)\left( \dfrac{\tilde{d}-\tilde{w}_{1}}{\tilde{w}_{1}}\right) \left( \dfrac{\tilde{d}-\tilde{w}_{2}}{\tilde{w}_{2}}\right) \left( \dfrac{\tilde{d}-\tilde{w}_{3}}{\tilde{w}_{3}}\right) \left( \dfrac{\tilde{d}-\tilde{w}_{4}}{\tilde{w}_{4}}\right) \\
  & = \left(\dfrac{2d-d}{d}\right)(a_{1}-1)(a_{2}-1)\left(\dfrac{2d-w_{3}}{w_{3}}\right)\left(\dfrac{2d-d}{d}\right) \\
  & = (a_{1}-1)(a_{2}-1)(2a_{3}-1).
\end{align*}
\hfill$\square$
\end{pf}

\begin{rmk}
Again from a similar argument as the one given in Remark 3.1, one can find for the weight vector $\tilde{\mathbf{w}}=(\tilde{w}_{0},\tilde{w}_{1},\tilde{w}_{2},\tilde{w}_{3},\tilde{w}_{4})$ obtained in Equation (31)  a polynomial of  Brieskorn-Pham type:
\begin{equation}    g^{T}_{BP}=z_{0}^{2}+z_{1}^{a_{1}}+z_{2}^{a_{2}}+z_{3}^{2a_{3}}+z_{4}^{2}\in H^0(\mathbb{P}({\bf\tilde{w}}), \mathcal{O}(\tilde{d}))
\end{equation}
with degree $\tilde{d}=\deg(g^T)=2d$ with $g$ given as (28) with no conditions on the exponents of $g.$
\end{rmk}

Next, we show a result similar to Theorem 4.6: we produce homotopy 7-spheres with empty extremal Sasaki-Reeb cone. These examples are obtained as a perturbation of a Berglund-Hübsch dual link of a polynomial of the form 
$g= z_{0}^{2}+z_{1}^{a_{1}}+z_{2}^{a_{2}}+z_{3}^{a_{3}}+z_{3}z_{4}^{2}.$

\begin{thm}   
    Consider the invertible polynomial $$g= z_{0}^{2}+z_{1}^{a_{1}}+z_{2}^{a_{2}}+z_{3}^{a_{3}}+z_{3}z_{4}^{2}$$ with $a_1\geq 4$ or $a_2\geq 4$ and $a_3\geq 3$. If the associated link $L_g$ satisfies the Lichnerowicz obstruction, then the Berglund-H\"ubsch dual 
    $L_{g^T}$ of $L_g$ satisfies this obstruction as well. Moreover, $L_{g^{T}_{BP}}$,  a perturbation of  $L_{g^T}$, which is defined via $$ g^{T}_{BP}=z_{0}^{2}+z_{1}^{a_{1}}+z_{2}^{a_{2}}+z_{3}^{2a_{3}}+z_{4}^{2}$$ given in Remark 4.4 satisfies the generalized Lichnerowicz obstruction (13) and therefore the link $L_{g_{BP}^{T}}$ does not admit extremal Sasaki metrics in its whole Sasaki-Reeb cone which is of dimension 2. Furthermore, if $a_0, a_1$ and $a_3$ are coprime and $a_{1},a_{2}$ are odd numbers, then $L_{g_{BP}^{T}}$ is a homotopy 7-sphere with no extremal metrics in its whole Sasaki-Reeb cone which is of dimension greater than one.
    \end{thm}

\begin{pf}
    For the polynomial $g$ defined in (28), we have its degree $d$ and its weight vector given by 
    \begin{equation}
        w_{0}=\dfrac{d}{2}, \ \ w_{1}=\dfrac{d}{a_{1}}, \ \ w_{2}=\dfrac{d}{a_{2}}, \ \ w_{3}=\dfrac{d}{a_{3}}, \ \ \mbox{ and } \ \ w_{4}=\dfrac{(a_{3}-1)d}{2a_{3}}
    \end{equation}
    where its index $I$ equals  
    $$I=|\textbf{w}|-d=w_{1}+w_{2}+\dfrac{w_{3}}{2}.$$
    On the other hand, from equality (31) we have 
{\small{
\begin{equation}
    \tilde{w}_{0}=d, \ \ \tilde{w}_{1}=2w_{1}=\dfrac{2d}{a_{1}}, \ \ \tilde{w}_{2}=2w_{2}=\dfrac{2d}{a_{2}}, \ \ \tilde{w}_{3}=w_{3}=\dfrac{d}{a_{3}}, \ \ \tilde{w}_{4}=d  
\end{equation}
}}
with  $\tilde{d}=2d.$

Notice that
$$\min_{i}\tilde{w}_{i}\in\{ \tilde{w}_{1},\tilde{w}_{2},\tilde{w}_{3}\}$$ and the index 
$\tilde{I}=|\tilde{\textbf{w}}|-\tilde{d}=\tilde{w}_{1}+\tilde{w}_{2}+\tilde{w}_{3}.$

Now, let us suppose that $L_{g}$ satisfies the Lichnerowicz obstruction. As $a_3\geq 3$ then ${w}_{4}\geq w_3.$ 
Thus $\min_{i}{w}_{i}\in\{{w}_{1},{w}_{2},{w}_{3}\}.$ 

\begin{itemize}
    \item[a)] Let us  suppose $\min_{i}{w}_{i}=w_1$ then 
        
        Since $L_{g}$ verifies the Lichnerowicz obstruction, we have
        $$I=w_{1}+w_{2}+\dfrac{w_{3}}{2}>4w_{1}.$$
         Multiplying by $2$ in the inequality above and the equalities given in (34), we obtain
        \begin{equation}         
        \tilde{I}=\tilde{w}_{1}+\tilde{w}_{2}+\tilde{w}_{3}=2w_{1}+2w_{2}+
        w_{3}>4(2w_{1})=4\tilde{w}_{1}\geq 4\min_{i}\tilde{w}_{i}  .          
        \end{equation}
         Thus, $L_{g^T}$ verifies the Lichnerowicz inequality. 
        \item[b)]  If $\min_{i}w_{i}=w_{2}$, the proof is similar to the argument given in case a).
         \item[c)]  Now, we assume that $\min_{i}w_{i}=w_{3}$. Since the weights $w_{i}$'s verify the Lichnerowicz inequality, we have
        \begin{equation*}
            I=w_{1}+w_{2}+\dfrac{w_{3}}{2}>4w_{3}.
        \end{equation*}
        From the equalities given (34) and the inequality given above, we obtain
        \begin{equation}
            \tilde{I}=\tilde{w}_{1}+\tilde{w}_{2}+\tilde{w}_{3}>w_{1}+w_{2}+\dfrac{w_{3}}{2}>4w_{3}=4\tilde{w}_{3}\geq 4\min_{i}\tilde{w}_{i}.
        \end{equation}
\end{itemize}
Thus, the link $L_{g^T}$ satisfies the Lichnerowicz inequality.  

It follows from the inequality (36) that 
         the  Brieskorn-Pham polynomial 
        $$g^{T}_{BP}=z_{0}^{2}+z_{1}^{a_{1}}+z_{2}^{a_{2}}+z_{3}^{2a_{3}}+z_{4}^{2}$$ defined in (32) satisfies the generalized Lichnerowicz  inequality (13). Therefore, 
        the link $L_{g_{BP}^{T}}$ does not admit extremal Sasaki metrics in its Sasaki-Reeb cone.
        The last claim follows from Lemma 4.7.
\hfill$\square$
 \end{pf}
    
\begin{rmk}
    Using a similar argument as the one given in Remark 4.3 we have the following: if we add an even number of quadratic terms $z_{i}^{2}$ to the polynomial $g_{BP}^{T}$:
    $$g_{Ext}^{T}=g_{BP}^{T}+z_{5}^{2}+\dots+z_{2n}^{2},$$
    then the link $L_{g_{Ext}^{T}}$ remains a homotopy sphere. Moreover, the dimension of associated cone increases by one for each two quadratic terms added, so we obtain a Sasaki-Reeb cone of dimension $n$. Also, the link does not admit  extremal Sasaki  metrics in its whole  Sasaki-Reeb cone.
\end{rmk}

As direct applications of Theorem 4.6 and Theorem 4.8 we give the following examples in dimension 7, which can be extrapolated to higher dimensions thanks to Remark 4.3 and Remark 4.5

\begin{example} We consider the polynomial $$f=z_{0}^{3}+z_{1}^{11}+z_{1}z_{2}^{2}+z_{3}^{13}+z_{3}z_{4}^{2}$$
 with weight vector $\mathbf{w}=(143,39,195,33,198)$, degree $d=429$ and index $I=179$. Applying the Berglund-Hübsch transpose rule, we obtain the polynomial 
$$ f^{T}  = z_{0}^{3}+z_{1}^{11}z_{2}+z_{2}^{2}+z_{3}^{13}z_{4}+z_{4}^{2}.$$ From this, we obtain 
 $$f_{BP}^{T}  = z_{0}^{3}+z_{1}^{22}+z_{2}^{2}+z_{3}^{26}+z_{4}^{2}.$$ Both polynomials 
associated to weight vector $\tilde{\mathbf{w}}=(286,39,429,33,429)$, with degree $\tilde{d}=858$ and index $\tilde{I}=358$. By Lemma 4.5, we have that the link $L_{f^{T}}$ is a rational homology sphere with $H_{3}(L_{f^{T}}, \mathbb{Z})=\mathbb{Z}_{3}$. Furthermore, we can verify that the  weight vector ${\bf w}$ satisfies  the Lichnerowicz inequality: $I=179>4(33)$. Therefore, the link $L_{f_{BP}^{T}}$ does not admit extremal Sasaki metrics.
\end{example}

\begin{example}
    We consider the polynomial  $$g=z_{0}^{2}+z_{1}^{5}+z_{2}^{7}+z_{3}^{12}+z_{3}z_{4}^{2}$$
with weight vector ${\bf w}=(420,168,120,70,385),$ degree $d=840$ and index $I=323$. If we apply the Berglund-Hübsch transpose rule, we obtain the polynomials
\begin{align*}
    g^{T} & =z_{0}^{2}+z_{1}^{5}+z_{2}^{7}+z_{3}^{12}z_{4}+z_{4}^{2}\\
    g_{BP}^{T} & = z_{0}^{2}+z_{1}^{5}+z_{2}^{7}+z_{3}^{24}+z_{4}^{2}
\end{align*}
with weights $\tilde{{\bf w}}=(420,168,120,35,420)$, degree $\tilde{d}=840$ and index $\tilde{I}=323$. By Lemma 4.7, we have that the link $L_{g_{BP}^{T}}$ is a homotopy 7-sphere. In addition, since the weight vector ${\bf w}$ satisfies  the Lichnerowicz inequality: $323>4(35)$, we get that the link $L_{g_{BP}^{T}}$ does not admit extremal Sasaki metrics.
\end{example}

\begin{example}
    For each $k\geq2$ integer, we define the polynomials
    $$f=z_{0}^{6k+1}+z_{1}^{2k-1}+z_{1}z_{2}^{2}+z_{3}^{2k+1}+z_{3}z_{4}^{2}$$
    We can verify that the integer numbers $6k+1, 2k-1$ and $2k+1$ are co-prime. The index $I$ is given by 
    $$I=w_{0}+\dfrac{w_{1}}{2}+\dfrac{w_{3}}{2}=\dfrac{d}{6k+1}+\dfrac{d}{2(2k-1)}+\dfrac{d}{2(2k+1)}$$
    Next, we will see that the weight vector ${\bf w}$ associated to $f$ verifies the Lichnerowicz inequality for all $k\geq2$.  Since $w_{0}=\min_{i}w_{i}$, we have 
    \begin{align*}
     I>4\min_{i}w_{i} 
 & \Longleftrightarrow \dfrac{d}{6k+1}+\dfrac{d}{2(2k-1)}+\dfrac{d}{2(2k+1)} >\dfrac{4d}{6k+1} \\
 & \Longleftrightarrow \dfrac{1}{2(2k-1)}+\dfrac{1}{2(2k+1)} >\dfrac{3}{6k+1}\\
 & \Longleftrightarrow \dfrac{4k}{4k^{2}-1} > \dfrac{6}{6k+1} \\
 & \Longleftrightarrow 24k^{2}+4k > 24k^{2}-6.
    \end{align*}
Thus, links $L_{f_{BP}^{T}}$ associated to polynomials of the form  $f_{BP}^{T}=z_{0}^{6k+1}+z_{1}^{2(2k-1)}+z_{2}^{2}+z_{3}^{2(2k+1)}+z_{4}^{2}$ are homotopy 
7-spheres and do not admit extremal Sasaki metrics.
\end{example}

\begin{example}
For each integer $k\geq2$, we consider the BP-chain polynomial:
$$g=z_{0}^{2}+z_{1}^{2k-1}+z_{2}^{2k+1}+z_{3}^{4k}+z_{3}z_{4}^{2}.$$
We can verify that the integers $2k-1, 2k+1$ and $4k$ are co-prime. Moreover, the index is given by
$$I=w_{1}+w_{2}+\dfrac{w_{3}}{2} = \dfrac{d}{2k-1}+\dfrac{d}{2k+1}+\dfrac{d}{8k}.$$
Now, we will show that the weights associated to  $g$ satisfy the Lichnerowicz inequality:
\begin{align*}
    I>4\min_{i}w_{i} & \Longleftrightarrow \dfrac{d}{2k-1}+\dfrac{d}{2k+1}+\dfrac{d}{8k} >4\left( \dfrac{d}{4k}\right) \\
    & \Longleftrightarrow \dfrac{1}{2k-1}+\dfrac{1}{2k+1} > \dfrac{7}{8k} \\
    & \Longleftrightarrow \dfrac{4k}{4k^{2}-1}>\dfrac{7}{8k}\\
    & \Longleftrightarrow 32k^{2}>28k^{2}-7.
\end{align*}
Thus, the links associated to the polynomials of the form $g_{BP}^{T}=z_{0}^{2}+z_{1}^{2k-1}+z_{2}^{2k+1}+z_{3}^{8k}+z_{4}^{2}$ are homotopy 7-spheres that do not admit extremal Sasaki metrics.
\end{example}

\begin{example}
    For distinct prime numbers $p,q$ and $r$ with $q\geq5$, $r\geq2$ and $p>rq$, define the polynomial
    $$f=z_{0}^{p^{k}}+z_{1}^{q^{k}}+z_{1}z_{2}^{2}+z_{3}^{r^{k}}+z_{3}z_{4}^{2},$$
    where $k\in\mathbb{Z}^{+}$. Then the link $L_{f_{BP}^T}$associated to the polynomial $f_{BP}^{T}$ obtained through the Berglund-Hübsch rule: $$f_{BP}^{T}=z_{0}^{p^{k}}+z_{1}^{2q^{k}}+z_{2}^{2}+z_{3}^{2r^{k}}+z_{4}^{2}$$
    is a rational homology 7-sphere with $H_{3}(L_{f_{BP}^{T}},\mathbb{Z})=\mathbb{Z}_{p^{k}}$ which does not admit extremal Sasaki metrics. Indeed, it is sufficient to verify the Lichnerowicz inequality for the weight vector $\bf w$ associated to $f$. Since $\min_{i}w_{i}=w_{0}=\dfrac{d}{p^{k}}$ and the index $I$ is given by 
    $$I=w_{0}+\dfrac{w_{1}}{2}+\dfrac{w_{3}}{2} = \dfrac{d}{p^{k}}+\dfrac{d}{2q^{k}}+\dfrac{d}{2r^{k}},$$
     we obtain 
    \begin{align*}
        I>4\min_{i}w_{i} & \Longleftrightarrow \dfrac{d}{p^{k}}+\dfrac{d}{2q^{k}}+\dfrac{d}{2r^{k}} >\dfrac{4d}{p^{k}} \\
        & \Longleftrightarrow \dfrac{1}{q^{k}}+\dfrac{1}{r^{k}} >\dfrac{6}{p^{k}} \\
        & \Longleftrightarrow r^{k}+q^{k}> 6\left(\dfrac{rq}{p}\right)^{k}.
    \end{align*}
    As $q\geq5$, $r\geq2$ and $p>rq$, we obtain $r^{k}+q^{k}\geq7$ and $1>\dfrac{rq}{p}$. Therefore, the inequality above is true. Notice that we can construct this rational homology 7-sphere for any prime number $p\geq11$.
\end{example}

\subsection{Examples derived from Thom-Sebastiani sums of more general invertible polynomials} 
If the topology of the link is not prescribed, it is not difficult to  produce 
more examples where the  Berglund-H\"ubsch transpose rule  maps links satisfying the Lichnerowicz inequality  to links preserving this feature and moreover that possess a representative in its local moduli which satisfies the generalized Lichnerowicz inequality (13). We do this for links of dimension 7 that come from Thom-Sebastiani sums of  blocks of invertible polynomials for the following type of singularities: 
\begin{itemize}
\item Thom-Sebastiani sum of singularities of the form BP+BP+chain.
\item Thom-Sebastiani sum of singularities of the form BP+cycle +chain.
\item Thom-Sebastiani sum of singularities of the form chain+chain.
\item Thom-Sebastiani sum of singularities of the form BP+chain.
\end{itemize}
The  proofs of the following lemmas are similar to the the ones given for Theorems 4.6 and 4.8 but we would rather include them for completeness. 

For the first case, we have the following result:
\begin{lemma}
    Given the invertible polynomial of the form
    $$f=z_{0}^2+z_{1}^{a_{1}}+z_{2}^{a_{2}}+z_{2}z_{3}^{a_{3}}+z_{3}z_{4}^2,$$
    where $a_{1}\geq 3$ and $a_{3}\geq 3$. If the corresponding  link $L_{f}$ satisfies the Lichnerowicz obstruction, then the Berglund-H\"ubsch dual $L_{f^T}$ also satisfies this obstruction. Furthermore, $L_{f^T}$ admits a perturbation $L_{f^{T}_{P}}$, associated to the polynomial $$f_{P}^{T}=z_{0}^2+z_{1}^{a_{1}}+z_{2}^{a_{2}}z_{3}+z_{3}^{2a_{3}}+z_{4}^2$$ with the same degree $\tilde{d}$ and same weight vector $\tilde{\bf{w}}$ obtained for $f^T$ such that it 
     does not admit extremal Sasaki metrics in its whole Sasaki-Reeb cone which is of dimension two.
\end{lemma}
    \begin{pf}
        Let $d$ be the degree of $f$. The weights $w_{i}$'s for the polynomial $f$ are given by
   {\small{
    \begin{equation}
        w_{0}=\dfrac{d}{2}, \ \\ w_{1}=\dfrac{d}{a_{1}}, \ \  w_{2}=\dfrac{d}{a_{2}} \ \  w_{3}=\dfrac{(a_{2}-1)d}{a_{2}a_{3}}, \ \ w_{4}=\dfrac{(a_{2}a_{3}-a_{2}+1)d}{2a_{2}a_{3}}
    \end{equation}
    }}
    where the index results
    $$I=\vert \textbf{w}\vert - d = w_{1}+w_{2}+\dfrac{w_{3}}{2}.$$
    The exponent matrix associated to $f$ is
    \begin{equation*}
       A_{f}= \begin{bmatrix}
          2 & 0 & 0 & 0 & 0  \\
      0 & a_{1} & 0 & 0 & 0 \\
      0 & 0 & a_{2} & 0 & 0 \\
      0 & 0 & 1 & a_{3} & 0 \\
      0 &  0 & 0 & 1 & 2   
        \end{bmatrix}
        \end{equation*}
    Applying the Berglund-Hübsch transpose rule, we obtain
    \begin{equation*}
       A^{T}_{f}= \begin{bmatrix}
          2 & 0 & 0 & 0 & 0  \\
      0 & a_{1} & 0 & 0 & 0 \\
      0 & 0 & a_{2} & 1 & 0 \\
      0 & 0 & 0 & a_{3} & 1 \\
      0 &  0 & 0 & 0 & 2   
        \end{bmatrix}
        \end{equation*}
    with associated polynomial $f^T=z_{0}^2+z_{1}^{a_{1}}+z_{2}^{a_{2}}z_{3}+z_{3}^{a_{3}}z_{4}+z_{4}^2$. The weight vector $\tilde{\textbf{w}}=(\tilde{w}_{1},\tilde{w}_{2},\tilde{w}_{3},\tilde{w}_{4})$ and degree $\tilde{d}$ corresponding to $f^T$ is obtained by solving the matrix equation
    \begin{equation*}
        \begin{bmatrix}
          2 & 0 & 0 & 0 & 0 & -1 \\
          0 & a_{1} & 0 & 0 & 0 & -1\\
          0 & 0 & a_{2} & 1 & 0 & -1\\
          0 & 0 & 0 & a_{3} & 1 & -1\\
          0 &  0 & 0 & 0 & 2 & -1  
        \end{bmatrix} \begin{bmatrix}
            \Tilde{w}_{0} \\
            \Tilde{w}_{1} \\
            \Tilde{w}_{2} \\
            \Tilde{w}_{3} \\
            \Tilde{w}_{4} \\
            \Tilde{d}
        \end{bmatrix} = \begin{bmatrix}
            0 \\ 
            0 \\
            0 \\
            0 \\
            0
        \end{bmatrix}
    \end{equation*}
    A solution for this equation is given by
    \begin{equation}
        \tilde{w}_{0}=\dfrac{\tilde{d}}{2}, \ \ \ \tilde{w}_{1}=\dfrac{\tilde{d}}{a_{1}}, \ \ \ \tilde{w}_{2}=\dfrac{(2a_{3}-1)\tilde{d}}{2a_{2}a_{3}}, \ \ \ \tilde{w}_{3} = \dfrac{\tilde{d}}{2a_{3}}, \ \ \ \tilde{w}_{4}=\dfrac{\tilde{d}}{2},    \end{equation}
    where $\tilde{d}=a_{2}a_{3}d$ and the index is 
    $$\tilde{I} = \vert \tilde{\bf w}\vert-\tilde{d}= \tilde{w}_{1}+\tilde{w}_{2}+\tilde{w}_{3}.$$
    Moreover, using the equalities in (38), we find a relation between $\tilde{I}$ and $I$:
  $$\tilde{I}=\tilde{w}_{1}+\tilde{w}_{2}+\tilde{w}_{3}=\dfrac{\tilde{d}}{a_{1}}+ \dfrac{(2a_{3}-1)\tilde{d}}{2a_{2}a_{3}}+\dfrac{\tilde{d}}{2a_{3}}=a_{2}a_{3}\left( w_{1}+w_{2}+\dfrac{w_{3}}{2}\right) = a_{2}a_{3}I$$
    Now, we suppose that $\bf w$ verifies the Lichnerowicz inequality. Since $a_{3}\geq 3$, we get $w_{4}\geq w_{3}$. Therefore, we have $\min_{i}w_{i} \in\{ w_{1},w_{2},w_{3}\}$. 
    
    Next, we will show that if $\bf w$ verifies the Lichnerowicz obstruction, its link $L_{f^T}$ also verifies this obstruction. 
    \begin{itemize}
        \item[a)]  If $\min_{i}w_{i}=w_{1}$, from our assumption we have the following 
        $$I=w_{1}+w_{2}+\dfrac{w_{3}}{2} > 4w_{1}.$$
        Multiplying by $a_{2}a_{3}$ the above inequality, we obtain
        $$\tilde{I}=a_{2}a_{3}I > 4a_{2}a_{3}w_{1} = \dfrac{4a_{2}a_{3}d}{a_{1}} = 4\tilde{w}_{1} \geq 4\min_{i}\tilde{w}_{i}$$
        Thus, the weight vector $\tilde{\bf w}$ verifies the Lichnerowicz inequality.
        \item[b)]  If $\min_{i}w_{i}=w_{2}$, we have
        $$I=w_{1}+w_{2}+\dfrac{w_{3}}{2} > 4w_{2} $$
        Multiplying by $a_{2}a_{3}$ again, we obtain
        $$\tilde{I}=a_{2}a_{3}I > 4a_{2}a_{3}w_{2} = 4a_{3}d$$
        In addition, from the equalities in (38), we have $$\tilde{w}_{2}=\dfrac{(2a_{3}-1)\tilde{d}}{2a_{2}a_{3}}=\dfrac{(2a_{3}-1)d}{2} < a_{3}d.$$ Replacing in the inequality above, we have 
        $$\tilde{I} > 4a_{3}d >4\tilde{w}_{2} \geq 4\min_{i}\tilde{w}_{i}.$$
        Thus, $\tilde{\textbf{w}}$ verifies Lichnerowicz obstruction.
        \item[c)]   If $\min_{i}w_{i}=w_{3}$, we have
        $$I=w_{1}+w_{2}+\dfrac{w_{3}}{2} > 4w_{3}. $$
        Once more, we multiply by $a_{2}a_{3}$ the last inequality. Then we obtain
        $$\tilde{I}=a_{2}a_{3}I > 4a_{2}a_{3}w_{3} = 4(a_{2}-1)d$$
        Also, from one of the equalities given in (38), we have $\tilde{w}_{3}=\dfrac{\tilde{d}}{2a_{3}}=\dfrac{a_{2}d}{2}$. As $a_{2}\geq2$, it verifies 
        $$\tilde{I}> 4(a_{2}-1)d \geq 4\left(\dfrac{a_{2}d}{2}\right) = 4\tilde{w}_{3} \geq 4\min_{i}\tilde{w}_{i}.$$
        Therefore, the weight vector $\tilde{\bf w}$ satisfies the Lichnerowicz inequality.
    \end{itemize}
    From these three cases we conclude that  $L_{f^T}$ preserves the Lichnerowicz obstruction. As  Remark 3.1,  $f^T$ can be perturbed to obtain a polynomial     $$f_{P}^{T}=z_{0}^2+z_{1}^{a_{1}}+z_{2}^{a_{2}}z_{3}+z_{3}^{2a_{3}}+z_{4}^2$$ with the same degree $\tilde{d}$ and same weight vector $\tilde{\bf{w}}$ obtained for $f^T$. Under this setting, it is clear  that the classical Lichnerowicz obstruction implies its generalized version on $L_{f^{T}_{P}}$. Thus the  link $L_{f^{T}_{P}}$ does not admit extremal Sasaki metrics in its whole Sasaki-Reeb cone.
 \hfill$\square$
 \end{pf}
\medskip

For  a Thom-Sebastiani sum of singularities of the form BP+cycle+chain we have the following: 

\begin{lemma}
    Given the invertible polynomial of the form
    $$f=z_{0}^2+z_{2}z_{1}^{a_{1}}+z_{1}z_{2}^{a_{2}}+z_{3}^{a_{3}}+z_{3}z_{4}^2,$$
    where $a_{3}\geq 3$. If the corresponding link $L_{f}$ sastisfies the Lichnerowicz obstruction given, then the  Berglund-H\"ubsch dual $L_{f^T}$ also satisfies this obstruction. Furthermore, $L_{f^T}$ admits a perturbation $L_{f^{T}_{P}}$ associated to the polynomial $$f_{P}^{T}=z_{0}^2+z_{2}z_{1}^{a_{1}}+z_{1}z_{2}^{a_{2}}+z_{3}^{2a_{3}}+z_{4}^2,$$ with the same degree $\tilde{d}$ and same weight vector $\tilde{\bf{w}}$ obtained for $f^T$, that does not admit extremal Sasaki metrics in its whole Sasaki-Reeb cone which is of dimension two.
    \end{lemma}
    \begin{pf}
        For the polynomial $f$, we have its associated weight vector $\textbf{w}=(w_{0},w_{1},w_{2},w_{3},w_{4})$ and degree $d$. Here $\bf w$ and $d$ verify
        {\small{
        \begin{equation}
            w_{0}= \dfrac{d}{2}, \ \ w_{1}=\dfrac{(a_{2}-1)d}{a_{1}a_{2}-1}, \ \ w_{2}=\dfrac{(a_{1}-1)d}{a_{1}a_{2}-1}, \ \ w_{3}=\dfrac{d}{a_{3}},  \ \ w_{4}=\dfrac{(a_{3}-1)d}{2a_{3}}
        \end{equation}
        }}
        where the index is
        $$I=\vert \textbf{w}\vert -d = w_{1}+w_{2}+\dfrac{w_{3}}{2}.$$
        For the polynomial $f$, its exponent matrix results
        \begin{equation*}
       A_{f}= \begin{bmatrix}
          2 & 0 & 0 & 0 & 0  \\
      0 & a_{1} & 1 & 0 & 0 \\
      0 & 1 & a_{2} & 0 & 0 \\
      0 & 0 & 0 & a_{3} & 0 \\
      0 &  0 & 0 & 1 & 2   
        \end{bmatrix}
        \end{equation*}
    Applying the Berglund-Hübsch transpose rule, we obtain
    \begin{equation*}
       A^{T}_{f}= \begin{bmatrix}
          2 & 0 & 0 & 0 & 0  \\
      0 & a_{1} & 1 & 0 & 0 \\
      0 & 1 & a_{2} & 1 & 0 \\
      0 & 0 & 0 & a_{3} & 1 \\
      0 &  0 & 0 & 0 & 2   
        \end{bmatrix}
        \end{equation*}
    with associated polynomial $f^T=z_{0}^2+z_{2}z_{1}^{a_{1}}+z_{1}z_{2}^{a_{2}}+z_{3}^{a_{3}}z_{4}+z_{4}^2$. In this case, the weight vector $\tilde{\bf w}$ associated to $f^T$ is given by
    $$\tilde{w}_{0}=\dfrac{\tilde{d}}{2}, \ \ \ \tilde{w}_{1}=\dfrac{(a_{2}-1)\tilde{d}}{a_{1}a_{2}-1}, \ \ \ \tilde{w}_{2} = \dfrac{(a_{1}-1)\tilde{d}}{a_{1}a_{2}-1}, \ \ \ \tilde{w}_{3}=\dfrac{\tilde{d}}{2a_{3}}, \ \ \ \tilde{w}_{4}=\dfrac{\tilde{d}}{2}$$
    where the degree is $\tilde{d}=d$. Thus, the index is 
    $$\tilde{I}=\vert \tilde{\bf w}\vert - \tilde{d}= \tilde{w}_{1}+\tilde{w}_{2}+\tilde{w}_{3}.$$
    Since $\tilde{w}_{1}=w_{1}$, $\tilde{w}_{2}=w_{2}$ and $\tilde{w}_{3}=\dfrac{w_{3}}{2}$, we have that $\tilde{I}=I$.
    Furthermore, as $a_{3}\geq 3$, then $\min_{i}w_{i}\in\{ w_{1},w_{2},w_{3}\}$.
    
    Now, we will prove that the Lichnerowicz inequality is preserved through the Berglund-Hübsch transpose rule. Indeed, we have 
    \begin{itemize}
        \item[a)] If $\min_{i}w_{i}=w_{1}$, the Lichnerowicz inequality is written as
        $$I=w_{1}+w_{2}+\dfrac{w_{3}}{2} > 4w_{1}.$$
        As $\tilde{I}=I$ and $\tilde{w}_{1}=w_{1}$, the inequality above implies
        $$\tilde{I} >4\tilde{w}_{1} \geq 4\min_{i}\tilde{w}_{i}.$$
        Therefore, $\tilde{\bf w}$ satisfies the Lichnerowicz obstruction.
        \item[b)] If $\min_{i}w_{i}=w_{2}$, the argument is similar to a). 
        \item[c)] If $\min_{i}w_{i}=w_{3}$, we have
        $$I=w_{1}+w_{2}+\dfrac{w_{3}}{2} >4w_{3}$$
        As $\tilde{I}=I$ and $w_{3}=2\tilde{w}_{3}$, we obtain
        $$\tilde{I} > 4w_{3} >4\tilde{w}_{3} \geq 4\min_{i}\tilde{w}_{i}.$$
        Thus, we have proved the Lichnerowicz inequality for $\tilde{\bf w}$.
    \end{itemize}
    From these three cases we conclude that  $L_{f^T}$ preserves the Lichnerowicz obstruction. 
    In addition, as in Remark 3.1,  we can always perturb  $f^T$ to  obtain  the  link $L_{f_{P}^T}$ associated to the polynomial:    $$f_{P}^{T}=z_{0}^2+z_{2}z_{1}^{a_{1}}+z_{1}z_{2}^{a_{2}}+z_{3}^{2a_{3}}+z_{4}^2,$$ with the same degree $\tilde{d}$ and same weight vector $\tilde{\bf{w}}$ obtained for $f^T$. It is obvious  that the classical Lichnerowicz obstruction implies its generalized version on this link. Thus the  link $L_{f^{T}_{P}}$ does not admit extremal Sasaki metrics in its whole Sasaki-Reeb cone.
\hfill$\square$
\end{pf}

For a Thom-Sebastiani sum of singularities of the form chain+chain we have the following: 
\begin{lemma}
    Let $f$ be the invertible polynomial    $$f=z_{0}^{a_{0}}+z_{0}z_{1}^{a_{1}}+z_{1}z_{2}^{2}+z_{3}^{a_{3}}+z_{3}z_{4}^2,$$
    where $a_{0},a_{1}$ and $a_{3}$ are greater than $2$. If the corresponding 
    link $L_{f}$ sastisfies the Lichnerowicz obstruction, then the Berglund-H\"ubsch dual $L_{f^T}$ also satisfies this obstruction. In addition, $L_{f^T}$ admits a perturbation $L_{f^{T}_{P}}$ associated to the polynomial $$f_{P}^T=z_{0}^{a_{0}}z_{1}+z_{1}^{2a_{1}}+z_{2}^{2}+z_{3}^{2a_{3}}+z_{4}^2$$
    with the same degree $\tilde{d}$ and same weight vector $\tilde{\bf{w}}$ obtained for $f^T$, that does not admit extremal Sasaki metrics in its whole Sasaki-Reeb cone which is of dimension two.
    \end{lemma}
        \begin{pf}
        If $d$ is the degree of polynomial $f$, then the weight vector $\bf w$ is given by
        $$w_{0}=\dfrac{d}{a_{0}}, \ \ \ w_{1}=\dfrac{(a_{0}-1)d}{a_{0}a_{1}}, \ \ \ w_{2}=\dfrac{(a_{0}a_{1}-a_{0}+1)d}{2a_{0}a_{1}}, \ \ \ w_{3}=\dfrac{d}{a_{3}}, \ \ \ w_{4}=\dfrac{(a_{3}-1)d}{2a_{3}},$$
        where the index is
        $$I=w_{0}+\dfrac{w_{1}}{2}+\dfrac{w_{3}}{2}.$$
        On the other hand, the exponent matrix of $f$ is 
        \begin{equation*}
       A_{f}= \begin{bmatrix}
          a_{0} & 0 & 0 & 0 & 0  \\
      1 & a_{1} & 0 & 0 & 0 \\
      0 & 1 & 2 & 0 & 0 \\
      0 & 0 & 0 & a_{3} & 0 \\
      0 &  0 & 0 & 1 & 2   
        \end{bmatrix}
        \end{equation*}
    Applying the Berglund-Hübsch transpose rule, we obtain
    \begin{equation*}
       A^{T}_{f}= \begin{bmatrix}
          a_{0} & 1 & 0 & 0 & 0  \\
      0 & a_{1} & 1 & 0 & 0 \\
      0 & 0 & 2 & 0 & 0 \\
      0 & 0 & 0 & a_{3} & 1 \\
      0 &  0 & 0 & 0 & 2   
        \end{bmatrix}
        \end{equation*}
    with associated polynomial $f^T=z_{0}^{a_{0}}z_{1}+z_{2}z_{1}^{a_{1}}+z_{2}^{2}+z_{3}^{a_{3}}z_{4}+z_{4}^2$. The weight vector $\tilde{\bf w}$ for $f^T$ is defined by
    $$\tilde{w}_{0}=\dfrac{(2a_{1}-1)\tilde{d}}{2a_{0}a_{1}}, \ \ \ \tilde{w}_{1}=\dfrac{\tilde{d}}{2a_{1}}, \ \ \ \tilde{w}_{2}=\dfrac{\tilde{d}}{2}, \ \ \ \tilde{w}_{3}=\dfrac{\tilde{d}}{2a_{3}}, \ \ \ \tilde{w}_{4}=\dfrac{\tilde{d}}{2},$$
    where $\tilde{d}=2a_{1}d$. Here, the index is 
    $$\tilde{I}=\vert \tilde{\bf w}\vert - \tilde{d} = \tilde{w}_{0}+\tilde{w}_{1}+\tilde{w}_{3}.$$
    Making a calculation, we obtain $\tilde{I}=2a_{1}I$. Furthermore, as $a_{1}\geq 3$ and $a_{3}\geq 3$, we have that $\min_{i}w_{i}\in \{ w_{0},w_{1},w_{3}\}$.

    Next, we will show that if we assume the Lichnerowicz inequality for $\bf w$, then $\tilde{\bf w}$ also verifies it. We consider the following cases:
    \begin{itemize}
        \item[a)] If $\min_{i}w_{i}=w_{0}$, then, from the assumption,  we have
        $$I=w_{0}+\dfrac{w_{1}}{2}+\dfrac{w_{3}}{2}>4w_{0}.$$
        Multiplying by $2a_{1}$ the previous inequality, we obtain
        $$\tilde{I}=2a_{1}I > 2a_{1}(4w_{0})= 4\left(\dfrac{2a_{1}d}{a_{0}}\right) > 4\left(\dfrac{(2a_{1}-1)d}{a_{0}}\right).$$
        Also, from  
        $$\tilde{w}_{0}=\dfrac{(2a_{1}-1)\tilde{d}}{2a_{0}a_{1}} = \dfrac{(2a_{1}-1)d}{a_{0}},$$
         we obtain $\tilde{I}>4\tilde{w}_{0}\geq 4\min_{i}\tilde{w}_{i}$. Thus, $\tilde{\bf w}$ verifies the Lichnerowicz obstruction.

        \item[b)] If $\min_{i}w_{i}=w_{1}$, then we have 
        $$I = w_{0}+\dfrac{w_{1}}{2}+\dfrac{w_{3}}{2} >4w_{1}.$$
        Multipying by $2a_{1}$ again, we obtain 
        $$\tilde{I}=2a_{1}I > 4(2a_{1}w_{1}) = 4\left(\dfrac{2(a_{0}-1)d}{a_{0}}\right) > 4d.$$
        Since $\tilde{w}_{1}=d$, we have $\tilde{I}>4\tilde{w}_{1}\geq 4\min_{i}\tilde{w}_{i}$. Therefore, the weight  vector $\tilde{\bf w}$ verifies the Lichnerowicz obstruction. 
        \item[c)] If $\min_{i}w_{i}=w_{3}$, then we have:
        $$I=w_{0}+\dfrac{w_{1}}{2}+\dfrac{w_{3}}{2} >4w_{3}.$$
        Multiplying by $2a_{1}$ we obtain
        $$\tilde{I}=2a_{1}I >4(2a_{1}w_{3}) = 4\left(\dfrac{2a_{1}d}{a_{3}}\right).$$
        Then, we notice that $\tilde{w}_{3}=\dfrac{\tilde{d}}{2a_{3}}=\dfrac{a_{1}d}{a_{3}}$. Thus $\tilde{I}>4\tilde{w}_{3}\geq 4\min_{i}\tilde{w}_{i}$.
    \end{itemize}
From these three cases we conclude that  $L_{f^T}$ preserves the Lichnerowicz obstruction. 
    The last statement of the lemma follows easily, arguing analogously as the previous lemma. 
    \hfill$\square$
 \end{pf}

For a Thom-Sebastiani sum of singularities of the form BP+chain we have the following: 

\begin{lemma}
    Consider an invertible polynomial of the form    $$f=z_{0}^2+z_{1}^{a_{1}}+z_{1}z_{2}^{a_{2}}+z_{2}z_{3}^{a_{3}}+z_{3}z_{4}^2,$$
    where $a_{1}\geq 2a_{3}\geq 6$. If the corresponding link $L_{f}$ sastisfies the Lichnerowicz obstruction, then the Berglund-H\"ubsch dual $L_{f^T}$ also satisfies this obstruction. In addition, 
    $L_{f^T}$ admits a perturbation $L_{f^{T}_{P}}$ associated to the polynomial 
    $$f^{T}_{P}=z_{0}^2+z_{1}^{a_{1}}z_{2}+z_{2}^{a_{2}}z_{3}+z_{3}^{2a_{3}}+z_{4}^2$$
    with the same degree $\tilde{d}$ and same weight vector $\tilde{\bf{w}}$ obtained for $f^T$, that does not admit extremal Sasaki metrics in its whole Sasaki-Reeb cone which is of dimension two.
\end{lemma}

    \begin{pf}
        The weight vector $\bf w$ associated to $f$ is defined by
        $$w_{0}=\dfrac{d}{2},\ \ \ w_{1}=\dfrac{d}{a_{1}}, \ \ \ w_{2}=\dfrac{(a_{1}-1)d}{a_{1}a_{2}},$$  $$w_{3}=\dfrac{(a_{1}a_{2}-a_{1}+1)d}{a_{1}a_{2}a_{3}}, \\\ w_{4}=\dfrac{(a_{1}a_{2}a_{3}-a_{1}a_{2}+a_{1}-1)d}{2a_{1}a_{2}a_{3}},$$
         where the index is
        $$I=w_{1}+w_{2}+\dfrac{w_{3}}{2}.$$
        Moreover, the exponential matrix associated to $f$ is
        \begin{equation*}
       A_{f}= \begin{bmatrix}
          2 & 0 & 0 & 0 & 0  \\
      0 & a_{1} & 0 & 0 & 0 \\
      0 & 1 & a_{2} & 0 & 0 \\
      0 & 0 & 1 & a_{3} & 0 \\
      0 &  0 & 0 & 1 & 2   
        \end{bmatrix}
        \end{equation*}
    Applying the Berglund-Hübsch transpose rule, we obtain
    \begin{equation*}
       A^{T}_{f}= \begin{bmatrix}
          2 & 0 & 0 & 0 & 0  \\
      0 & a_{1} & 1 & 0 & 0 \\
      0 & 0 & a_{2} & 1 & 0 \\
      0 & 0 & 0 & a_{3} & 1 \\
      0 &  0 & 0 & 0 & 2   
        \end{bmatrix}
        \end{equation*}
    with associated polynomial $f^T=z_{0}^2+z_{1}^{a_{1}}z_{2}+z_{2}^{a_{2}}z_{3}+z_{3}^{a_{3}}z_{4}+z_{4}^2$. Here, the weight vector $\tilde{\bf w}$ is given by
    $$\tilde{w}_{0}=\dfrac{\tilde{d}}{2}, \ \ \ \tilde{w}_{1}=\dfrac{(2a_{2}a_{3}-2a_{3}+1)\tilde{d}}{2a_{1}a_{2}a_{3}}, \ \ \ \tilde{w}_{2}=\dfrac{(2a_{3}-1)\tilde{d}}{2a_{2}a_{3}}, \ \ \ \tilde{w}_{3}=\dfrac{\tilde{d}}{2a_{3}}, \ \ \ \tilde{w}_{4}=\dfrac{\tilde{d}}{2},$$
    where $\tilde{d}=a_{1}a_{2}a_{3}d$ and the index $\tilde{I}$ is given by 
$$\tilde{I}=\vert \tilde{\bf w}\vert -\tilde{d} = \tilde{w}_{1}+\tilde{w}_{2}+\tilde{w}_{3}.$$
Some calculations lead to  $\tilde{I}=a_{1}a_{2}a_{3}I$. Moreover, since $a_{3}>2$, we have  $w_{3}\leq w_{4}$. Thus,   $\min_{i}w_{i}\in\{ w_{1},w_{2},w_{3}\}$.

\begin{itemize}
    \item[a)] If $\min_{i}w_{i}=w_{1}$, the Lichnerowicz inequality is given by 
    $$I=w_{1}+w_{2}+\dfrac{w_{3}}{2} >4w_{1}$$
    Multiplying by $a_{1}a_{2}a_{3}$, we obtain 
    $$\tilde{I}=a_{1}a_{2}a_{3}I>4a_{1}a_{2}a_{3}w_{1}=4a_{2}a_{3}d.$$
    On the other hand, as $\tilde{d}=a_{1}a_{2}a_{3}d$, we have
    $$\tilde{w}_{1}=\dfrac{(2a_{2}a_{3}-2a_{3}+1)\tilde{d}}{2a_{1}a_{2}a_{3}}=\dfrac{(2a_{2}a_{3}-2a_{3}+1)d}{2}<a_{2}a_{3}d$$
    Replacing above, we get $\tilde{I}>4\tilde{w}_{1}\geq \min_{i}\tilde{w}_{i}$.
\item[b)] If $\min_{i}w_{i}=w_{2}$, then 
    $$I=w_{1}+w_{2}+\dfrac{w_{3}}{2}>4w_{2}.$$ 
    Multiplying by $a_{1}a_{2}a_{3}$, we obtain 
    $$\tilde{I}=a_{1}a_{2}a_{3}I > 4a_{1}a_{2}a_{3}w_{2} = 4a_{3}(a_{1}-1)d.$$
    Also, as $a_{1}\geq 2a_{3}$, we have $a_{3}(a_{1}-1) \geq \dfrac{a_{1}(2a_{3}-1)}{2}$. Replacing above, we obtain
    $$\tilde{I} > 4 \left(\dfrac{a_{1}(2a_{3}-1)d}{2}\right) = 4\left(\dfrac{(2a_{3}-1)\tilde{d}}{2a_{2}a_{3}}\right) = 4\tilde{w}_{2} \geq 4\min_{i}\tilde{w}_{i}.$$
    Thus  the weight vector $\tilde{\bf w}$ verifies the Lichnerowicz inequality.

    \item[c)] If $\min_{i}w_{i}=w_{3}$, we have
    $$I= w_{1}+w_{2}+\dfrac{w_{3}}{2}>4w_{3}.$$
    Once more, multiplying by $a_{1}a_{2}a_{3}$ above, we obtain
    $$\tilde{I}=a_{1}a_{2}a_{3}I > 4a_{1}a_{2}a_{3}w_{3} = 4(a_{1}a_{2}-a_{1}+1)d.$$
    On the other hand, as $\tilde{w}_{3}=\dfrac{\tilde{d}}{2a_{3}}=\dfrac{a_{1}a_{2}d}{2}$, 
    $$\tilde{I}>4(a_{1}a_{2}-a_{1}+1)d>4\left(\dfrac{a_{1}a_{2}d}{2}\right) = 4\tilde{w}_{3}\geq 4\min_{i}w_{i}.$$
    Thus, $\tilde{\bf w}$ satisfies the Lichnerowicz inequality.
\end{itemize}

From these three cases we conclude that  $L_{f^T}$ preserves the Lichnerowicz obstruction. 
    The last statement of the theorem follows from similar arguments as the ones used in the previous lemmas. 
   \hfill$\square$
    \end{pf}

\section*{Declarations}





\subsection*{Funding} The first author received  financial support from Pontificia Universidad Católica del Perú through project DGI 2019-1-0089. 

\subsection*{Data Availability} Data sharing not applicable to this article as no datasets were generated or analyzed during the current study.

\subsection*{Conflict of interests} The authors have no competing interests to declare that are relevant to the content of this article.



\end{document}